\documentclass[11pt]{article} 

\usepackage[utf8]{inputenc} 

\usepackage{amssymb,amsmath,mathtools}
\usepackage{graphicx}
\usepackage{geometry}
\usepackage{subcaption}
\usepackage{array}
\usepackage{makecell,multirow}
\usepackage{authblk}
\newcommand{\R}{\mathbb{R}}
\newcommand{\Rad}{\mathcal{R}}

\newcommand{\Radh}{\mathcal{R}_{\text{h}}}
\newcommand{\Oc}{\mathcal{O}}
\newcommand{\betah}{\frac{\beta}{2}}
\newcommand{\Fc}{\mathcal{F}}
\graphicspath{{./fig/}}

\def\Plus{\texttt{+}}
\def\Minus{\texttt{-}}
\mathtoolsset{showonlyrefs=true}

\date{}
\title{Fast hyperbolic Radon transform represented as convolutions in log-polar coordinates}
\author[1]{Viktor~V.~Nikitin}

\author[1]{Fredrik~Andersson}

\author[1]{Marcus~Carlsson}

\author[2]{Anton~A.~Duchkov}
\affil[1]{Centre for Mathematical Sciences, Lund University,
	S\"olvegatan 18, Box 118, SE-22100 Lund,~Sweden}
\affil[2]{Institute of Petroleum Geology and Geophysics SB RAS, 
3, Ac. Koptyuga ave., 630090 Novosibirsk,~Russian~Federation}

\begin{document}
	\maketitle			
\begin{abstract}
	The hyperbolic Radon transform is a commonly used tool in seismic processing, for instance in seismic velocity analysis, data interpolation and for multiple removal. A direct implementation by summation of traces with different moveouts is computationally expensive for large data sets. In this paper we present a new method for fast computation of the hyperbolic Radon transforms. It is based on using a log-polar sampling with which the main computational parts reduce to computing convolutions. This allows for fast implementations by means of FFT. In addition to the FFT operations, interpolation procedures are required for switching between coordinates in the time-offset;  Radon; and log-polar domains. Graphical Processor Units (GPUs) are suitable to use as a computational platform for this purpose, due to the hardware supported interpolation routines as well as optimized routines for FFT.
	Performance tests show large speed-ups of the proposed algorithm. Hence, it is suitable to use in iterative methods, and we provide examples for data interpolation and multiple removal using this approach.
\end{abstract}



\section{Introduction}
In the processing of Common-Midpoint gathers (CMPs), the hyperbolic Radon transform has proven to be a valuable tool for instance in velocity analysis \cite{clayton1981inversion,greenhalgh1990seismic}; aliasing and noise removal \cite{turner1990aliasing}; trace interpolation \cite{averbuch2001fast,yu2007wavelet}; and attenuation of multiple reflections \cite{hampson1986inverse}. The hyperbolic Radon transform is defined as
\begin{equation}
\Radh f(\tau,q)=\int_{-\infty}^{\infty}f\left(\sqrt{\tau^2+q^2 x^2},x\right)d x,\label{Rhdef}
\end{equation}
where the function $f(t,x)$ usually corresponds to a CMP gather. Here, the parameter $q$ characterizes an effective velocity value; and $\tau$ represents the intercept time at zero offset.

Several versions of Radon transforms are used in seismic processing, e.g., straight-line, parabolic, and hyperbolic Radon transforms. In many applications there is a need for a sparse representation of seismic data using hyperbolic wave events. One way to get sparse representations is by using iterative thresholding algorithms with sparsity constraints \cite{daubechies2004iterative,sacchi1995high}. Popular applications using such representations are seismic data interpolation and wavefield separation \cite{jiang2016time,trad2003interpolation}.
Since iterative schemes for computing such representations require the application of the forward and adjoint operators several times, it becomes important to use fast algorithms for the realization of the operators to the limit the total computational cost.

Note that the direct summation over hyperbolas in \eqref{Rhdef} has a computational complexity of $\Oc(N^3)$, given that the numbers of samples for the variables $t,x,\tau,q$ are $\mathcal{O}(N)$. There are many effective ($\Oc(N^2\log N)$) methods for rapid evaluation of the traditional Radon transforms, or the parabolic Radon transform, see \cite{beylkin1984inversion,fessler2003nonuniform,schonewille2001parabolic}.
The hyperbolic Radon transform is, however, more challenging. Nonetheless, a fast method for hyperbolic Radon transforms was recently presented in \cite{hu2013fast}. This method is based on using the fast butterfly algorithms described in \cite{candes2009fast,o2007new,poulson2014parallel}.

A fast method for the standard Radon transform was proposed in \cite{andersson2005fast} by expressing the Radon transform and its adjoint in terms of convolutions in log-polar coordinates. Computationally efficient algorithms for GPUs were presented for this approach in \cite{andersson2015fast}. In this paper we propose to use the same approach and construct algorithms with complexity $\Oc(N^2 \log N)$ for evaluation of the hyperbolic Radon transform. We present computational performance tests confirming the expected accuracy and the computational complexity, as well as predicted computational speed-ups for parallel implementations. Finally, we present several synthetic and real data tests using the hyperbolic Radon transform for data interpolation and multiple attenuation.

\section{Method}
To begin with, we note that functions $f(t,x)$ describing CMP gathers are symmetric with respect to $x=0$. Hence, by introducing \begin{equation}\label{f}\tilde{f}(s,y)=\frac{f(\sqrt{s},\sqrt{y})}{2\sqrt{y}},\end{equation} it follows that
\begin{equation}\label{Rhdef2}\begin{aligned}
\Radh f(\tau,q)=
2\int_{0}^{\infty}f\left(\sqrt{\tau^2+q^2 x^2},x\right)dx=2\int_{0}^{\infty}\tilde{f}(\tau^2+q^2 y,y)dy.\end{aligned}
\end{equation}
The resulting expression in \eqref{Rhdef2} has a form of the Radon transform over straight lines, and a fast algorithm for the evaluation of this was presented in \cite{andersson2015fast}, referred to as the log-polar Radon transform which is based on rewriting the key operations as convolutions in a log-polar coordinate system. In Section \ref{logpol} we briefly recall the construction of the log-polar Radon transform and discuss how to adjust this method for optimal performance when processing seismic data, and in Section \ref{hyp} we introduce coordinate transforms as well as sampling/interpolation requirements for accurate evaluation of $\Radh f (\tau q)$.
\begin{figure}
	\centering 	 \includegraphics[width=0.425\textwidth]{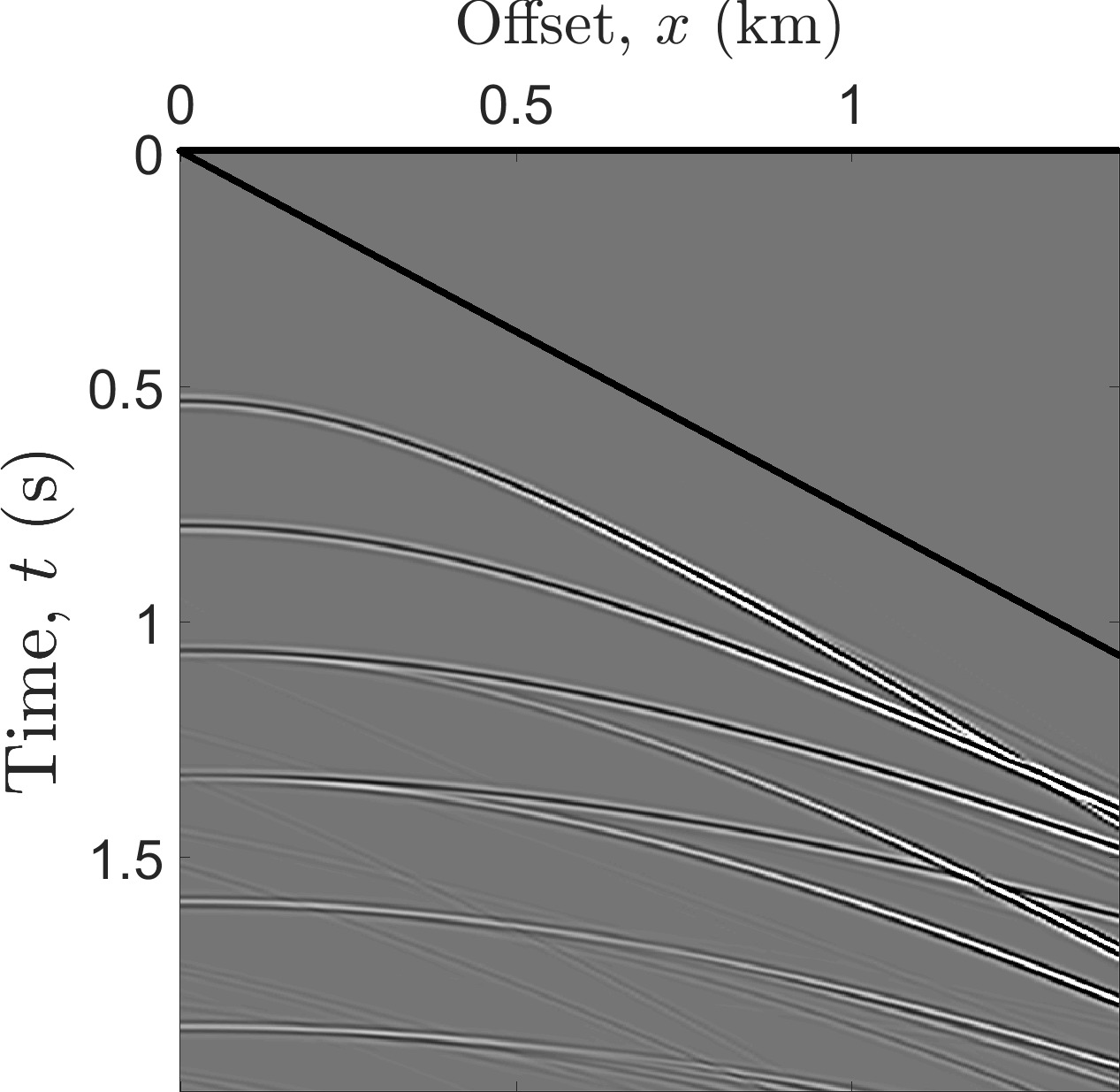}	\hspace{5pt}
	\includegraphics[width=0.425\textwidth]{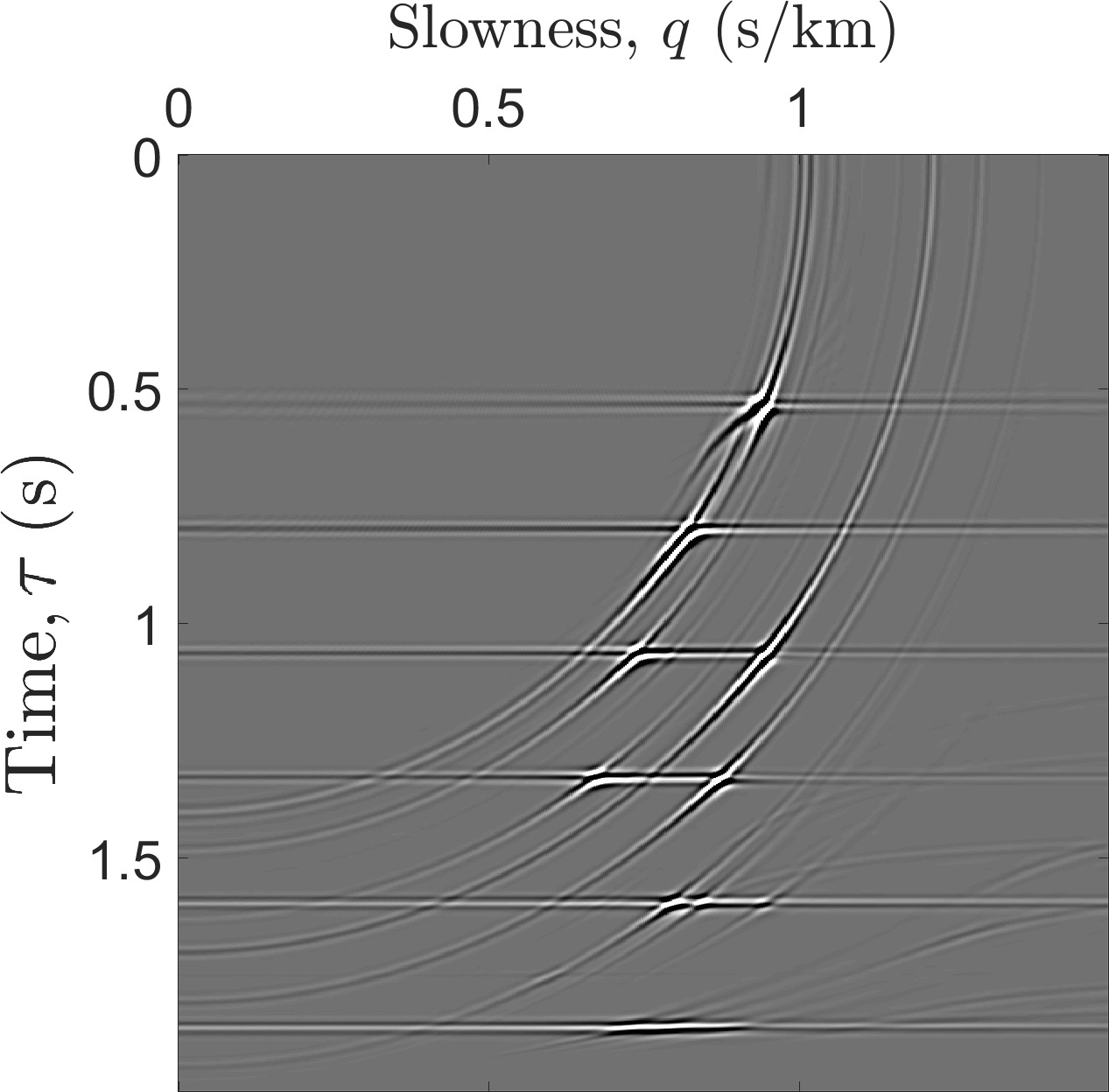}
	\caption{CMP gather, cutting information for high offset and small time intercept (left). Hyperbolic Radon transform, the circles are artifacts from truncation (right).}
	\label{fig:example}
\end{figure}

\subsection{Log-polar Radon transform}\label{logpol}
The standard Radon transform (cf. \eqref{Rhdef2}) can be written in terms of a double integral
\begin{equation}\label{taupint}\begin{aligned}
\Rad \tilde{f}(\tau^2,q^2)=\iint \tilde{f}(s,y)\delta(s-\tau^2-q^2 y)dyds,
\end{aligned}
\end{equation}
where $\delta$ denotes the Dirac distribution. In \cite{andersson2015fast} one works with the log-polar coordinates
\begin{align}
\begin{cases}
s=e^{\rho'}\cos(\theta'),\\
y=e^{\rho'}\sin(\theta'),
\end{cases}\quad\quad
\begin{cases}
\tau^2=\frac{e^\rho}{\cos(\theta)},\\
q^2=-\tan{\theta}.
\end{cases}
\label{param}
\end{align}
By introducing $\zeta(\theta,\rho)=\delta(\cos(\theta)-e^{\rho})$, it turns out that the Radon transform can be efficiently evaluated using the \textit{log-polar Radon transform}
\begin{equation}\begin{aligned}\label{lpfwd}
\Rad_{\text{lp}}\tilde{f}(\theta,\rho)=
\cos(\theta)\iint \tilde{f}(\theta',\rho')e^{\rho'}\zeta(\theta-\theta',\rho-\rho')d\rho'd\theta'\end{aligned},
\end{equation}
where, by abuse of notation, we use the same notation $\tilde{f}$ for both coordinate representations.

\begin{figure}
	\centering 	
	\begin{subfigure}[b]{0.35\textwidth} {\includegraphics[width=1\textwidth]{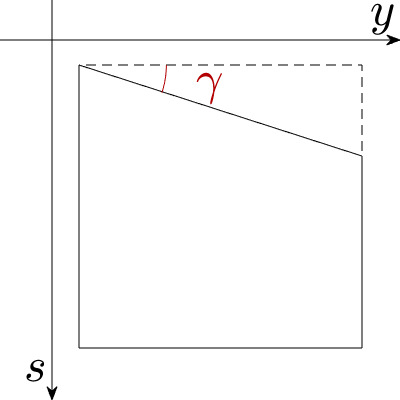}}	\end{subfigure}	\hspace{0.5cm}
	\begin{subfigure}[b]{0.53\textwidth}		 {\includegraphics[width=1\textwidth]{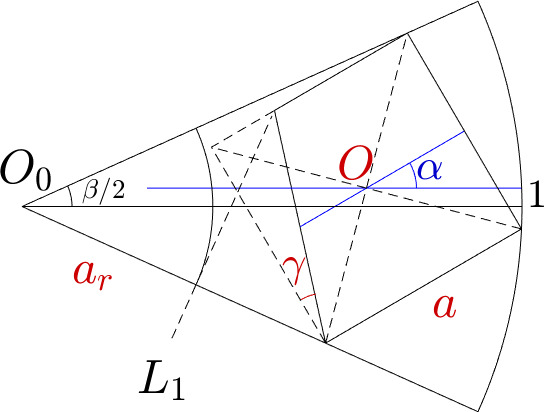}}\end{subfigure}	
	\caption{Region of interest (trapezoid) for data evaluation (left). Scaling, rotation and translation for the log-polar setup (right).}
	\label{fig:setupRp}
\end{figure}

However, the above representation is not suitable for treating functions $\tilde{f}$ with support near 0, since this corresponds to $\rho'=-\infty$. One therefore applies scaling, rotation and translation to work with functions supported within a subset of a circle-sector of opening angle $\beta$ as in Figure \ref{fig:setupRp}, right. Due to certain technicalities \cite{andersson2015fast}, the implementation of $\Rad_{\text{lp}}$ works best when evaluating only for values $\theta\in [-\beta/2,\beta/2]$. We will refer to this algorithm as the \textit{partial} $\Rad_{\text{lp}}$.

With this in mind, we now briefly explain how to make slight modifications to the above scheme, better suited for the processing of CMP gathers. A simplified synthetic example of a typical CMP gather is shown in Figure \ref{fig:example}. Note that the function continues outside the maximum limits given by $x$ and $t$, leading to a truncation of \eqref{Rhdef2}, (which can be seen e.g. as the circular artifacts in Figure \ref{fig:example}). Also note that there is no data in the region above a line $t=kx$, i.e. high offset $x$ and small time intercept $t$, so to decrease the amount of computations we may ignore this piece. In the coordinates $(s,y)$ this triangle is again a triangle, but with equation $s=k^2y$. We set $\gamma=\arctan k^2$. Thus, we are in practice only interested in evaluating \eqref{taupint} for data $\tilde{f}$ on a right trapezoid with the form illustrated in Figure \ref{fig:setupRp}, left. Besides, one is usually also interested only in values of $(\tau,q)$ in a limited range $[\tau_{min},\tau_{max}]\times [q_{min},q_{max}]$.

In order for this to correspond to a symmetric interval of $\theta$, we set $\beta=\arctan(q^2_{\max})-\arctan(q^2_{\min})$ and modify the relation between $\theta$ and $q$ in \eqref{param} as follows $$\theta=\alpha-\arctan(q^2),$$ where
$\alpha=(\arctan(q^2_{\max})+\arctan(q^2_{\min}))/2$. For a particular value of $\theta$ the output of the partial $\Rad_{\text{lp}}$  correspond to integrals over lines whose angle with respect to the vertical axis is $\theta$. In order for these to correspond to desired values of $q$, one needs to rotate $\tilde{f}$ so that the $t$ axis makes an angle $\alpha$ with respect to the horizontal axis in Figure \ref{fig:setupRp}, right. Moreover, due to the problems at the origin, $\tilde{f}$ needs to be dilated and translated so that it fits within the circle sector of radius 1 and opening angle $\beta$, as in Figure \ref{fig:setupRp}, right. This has the effect that the trapezoidal support is inscribed inside a square with side length $a$, located so that three of its corners lie on the border of the sector. It can be shown that
\begin{align*}
& a=\frac{\sin (\beta)}{\sqrt{\sin (2 \alpha ) \sin (\beta )+\cos (\beta ) (\sin (2 \alpha )+\sin (\beta ))+1}},\\
& O=O\left(\frac{a \sin \left(\alpha +\frac{\pi }{4}\right) \tan (\beta )}{\sqrt{2}},\frac{a \cos \left(\alpha +\frac{\pi }{4}\right) \tan \left(\frac{\beta }{2}\right)}{\sqrt{2}}\right).
\end{align*}
The line $L_1$ passes through the fourth corner of the trapezoid and is orthogonal to the border of the sector; the distance from the origin $O_0$ to the line $L_1$ is indicated by $a_r$ and indicates the the first non-zero contribution to the partial $\Rad_{\text{lp}}$.

\begin{figure}
	\setlength{\extrarowheight}{52pt}
	\centering 	
	\begin{tabular}{ >{\centering\arraybackslash}m{1in} c >{\centering\arraybackslash}m{1in} c}
		\begin{subfigure}[b]{0.17\textwidth}{	\centering 	\includegraphics[width=1\textwidth]{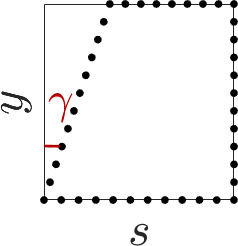}}\end{subfigure}& {\LARGE$\xrightarrow{}$}&
		\begin{subfigure}[b]{0.17\textwidth}{	\centering 	\includegraphics[width=1\textwidth]{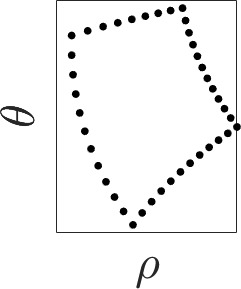}}\end{subfigure} &
		\multirow{-1}[4]{*}{	 \begin{subfigure}[b]{0.4\textwidth}{\includegraphics[width=1\textwidth]{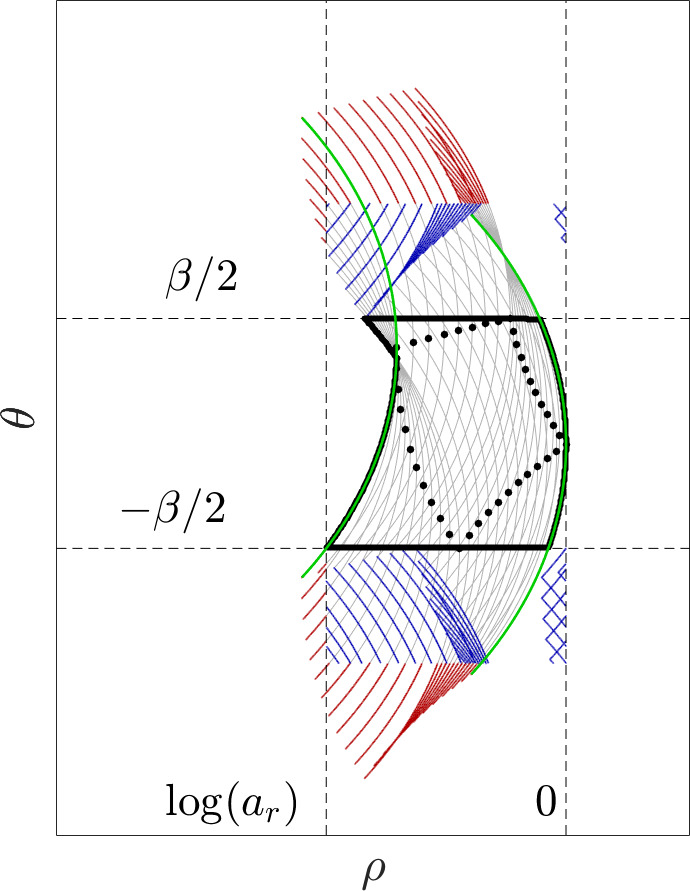}}\hspace{1cm}\end{subfigure}}\\
		
		\begin{subfigure}[b]{0.17\textwidth}{	\centering 	\includegraphics[width=1\textwidth]{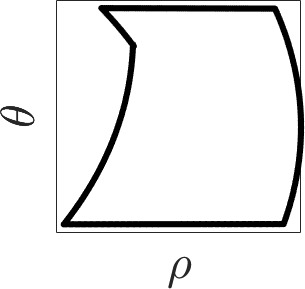}}\end{subfigure}&
		{\LARGE$\xrightarrow{}$}&
		\begin{subfigure}[b]{0.17\textwidth}{	\centering 	\includegraphics[width=1\textwidth]{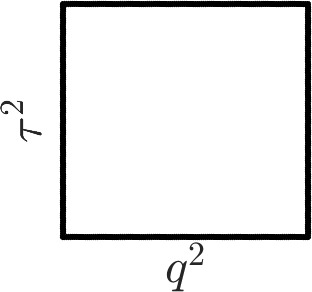}}\end{subfigure}&
	\end{tabular}
	\caption{Conversion to the log-polar domain and effects of computing convolutions.}
	\label{fig:convs}
\end{figure}

In summary, we are interested in the values of the log-polar Radon transform in the range $\left[-\frac{\beta}{2},\frac{\beta}{2}\right]\times[\log(a_r),0]$. With this setup the log-polar Radon transform can be computed in terms of the finite convolution
\begin{equation}\label{finiteconv}\begin{aligned}
\Rad_\text{lp}\tilde{f}(\theta,\rho)=\cos(\theta)\int_{-\frac{\beta}{2}}^{\frac{\beta}{2}}\int_{\log(a_r)}^{0}&\tilde{f}(\theta',\rho')e^{\rho'}\zeta(\theta-\theta',\rho-\rho')d\rho'd\theta'=\\&\cos(\theta)\Fc^{\text{-}1}\Big(\Fc \left(\tilde{f}(\theta,\rho)e^\rho\right)(\hat{\theta},\hat{\rho})\cdot\Fc \zeta(\hat{\theta},\hat{\rho})\Big)(\theta,\rho). \end{aligned}
\end{equation}
Here, $\Fc$ denotes the two-dimensional Fourier transform. We use the notation  $(\hat{\theta},\hat{\rho})$ for the reciprocal variables of $(\theta,\rho)$. The function $\widehat{\zeta}(\hat{\theta},\hat{\rho})$ can be accurately evaluated numerically (in a precomputing step) in contrast to $\zeta(\theta,\rho)$ which is defined in terms of distributions and is discontinuous along a curve, see \cite{andersson2015fast} for a detailed description. To avoid wrapping effects, zero-padding is applied in the log-polar domain. The effects of the convolutions are schematically  illustrated in Figure \ref{fig:convs}. The trapezoid containing the support of the data is transformed to the shape indicated by the black points after a change to log-polar coordinates; the green lines show shifted versions of the function $\zeta$; and the support after the log-polar Radon transform is applied is indicated by the thick black curves. By using this scheme we conclude that the rectangle $[-\beta,\beta]\times[\log(a_r),0]$ is a good choice for enclosing the support of the functions, which is needed for the discrete evaluation of the integrals by means of convolutions in log-polar coordinates.

We now describe how $\Rad_{\text{lp}}$ can be used to recover $\Rad \tilde{f}$ for a function $\tilde{f}$ with support in the unit rectangle.
The change of coordinates $(s,y)$ for the log-polar setup is described by the transformation $\mathsf{T}$,
\begin{align}\label{transformT}
\mathsf{T}\begin{pmatrix}t\\x\end{pmatrix}=a\begin{pmatrix}\cos(\alpha)&-\sin(\alpha)\\ \sin(\alpha) &\cos(\alpha)\end{pmatrix}\begin{pmatrix}s-0.5\\y-0.5\end{pmatrix}+\begin{pmatrix}O_1\\O_2\end{pmatrix}
\end{align}
as well as the change of coordinates $(\tau,q)$ for the log-polar setup can be expressed by $\mathsf{S}$, which can be found by scaling, rotation and translation procedures for the log-polar setup. Some tedious manipulations yield
\begin{equation}\label{transformS}\begin{aligned}
\mathsf{S}\begin{pmatrix}\tau^2\\q^2
\end{pmatrix}=\begin{pmatrix}
\left(a(\tau^2-
\frac{1}{2})\cos(\alpha)+a\frac{\sin(\alpha)}{2}+O_1+\left(a(\tau^2-\frac{1}{2})\sin(\alpha)-a\frac{\cos(\alpha)}{2}+O_2\right)\!\phi\right)\\
\phi
\end{pmatrix}
\\\text{with } \phi=\tan(\alpha-\arctan(q^2)).
\end{aligned}
\end{equation}
Moreover, we introduce two transformations for switching to log-polar coordinates according to relations (\ref{param}):
\begin{align}\label{transformP}
&\mathsf{P}_1\begin{pmatrix}t\\x
\end{pmatrix}=
\begin{pmatrix}
\log(\sqrt{t^2+x^2})\\\arctan\left(\frac{x}{t}\right)\end{pmatrix}
&\mathsf{P}_2\begin{pmatrix}\tau^2\\q^2
\end{pmatrix}=
\begin{pmatrix}
\log(\tau^2\cos(-\arctan(q^2)))\\-\arctan\left(q^2\right)\end{pmatrix}
\end{align}
To the end, by introducing linear operators
\begin{equation}\label{coord_trans}
\begin{aligned}
T\tilde f=\tilde{f}\left(\mathsf{T}^{-1}\mathsf{P}_1^{-1}\cdot\right) \quad Sg=g\left(\mathsf{S}^{-1}\mathsf{P}_2^{-1}\cdot\right)
\end{aligned}
\end{equation}
the Radon transform over straight lines and its adjoint operator can be computed (up to a scaling factor) by
\begin{eqnarray}\label{Radon_all_fwd_adj}
\Rad \tilde{f}(\tau,q) =S^{-1} \Rad_{\text{lp}} \left(T \tilde{f}\right)\left(\tau,q\right),\\
\Rad^* g(t,x) =T^{-1}\Rad_{\text{lp}}^* \left(S g\right)\left(t,x\right).
\end{eqnarray}
\subsection{Hyperbolic coordinates}\label{hyp}
Let $f$ be a CMP gather measured on the rectangle
\begin{equation}\label{Rgrid}
\{(t,x): 0\le t\le T, 0\le x \le X\}.
\end{equation}
which we treat as a function on all of $\R^2$ which is 0 outside this rectangle.
Note that
\begin{equation}\begin{aligned}\label{hypnorm}
&\Radh (f)\left(\tau,q\right)=XT\int_{0}^{1}\int_{0}^{1} {f}(Tt,Xx)\delta\left(tT-\sqrt{\tau^2+q^2x^2 X^2}\right)d x d t=\\
&X\int_{0}^{1}\int_{0}^{1} {f}(Tt,Xx)\delta\left(t-\sqrt{(\tau/T)^2+q^2x^2 \frac{X^2}{T^2}}\right)d x d t=X\Radh \Big({f}(T\cdot,X\cdot)\Big)\left(\frac{\tau}{T},\frac{qX}{T}\right),
\end{aligned}
\end{equation}
which allows us to assume that $f$ is given on the rectangle $[0,1]\times[0,1]$ to begin with. Upon corresponding rescaling of $\tau:=\frac{\tau}{T}$ and $q:=\frac{qX}{T}$, we are interested in evaluating $\Radh f$ on the rectangle
\begin{equation}\label{Qgrid}
\{(\tau,q): \tau_{min}\le\tau\le 1, q_{\text{min}}\le q\le q_{\text{max}}\},
\end{equation}
where $\tau_{min}$ corresponds to the arrival of the first event in the rescaled coordinates.

Now we recall the expression \eqref{finiteconv} representing the finite convolution for computing the log-polar Radon transform. It can be rapidly evaluated in terms of FFT if the log-polar samples $(\theta,\rho)$ are given on an equally spaced grid. Since data is assumed to be sampled in the $(t,x)$ domain, a resampling is needed.
We propose to do this using  \textit{cardinal B-spline interpolation} \cite{deboor1978practical,unser1999splines}, since this type of interpolation is particularly well suited for GPU implementations (cf. \cite{ruijters2008efficient}). This technique is related to that used for fast unequally-spaced Fourier transforms (USFFT) \cite{USFFT,Rokhlin_Dutt}, in the way that the interpolation is conducted by smearing data in one of the domains, and the compensating for that effect is done in the reciprocal domain.

In \eqref{finiteconv} we have to compute $\Fc \left(\tilde f(\theta,\rho)e^\rho\right)(\hat{\theta},\hat{\rho})$ which we can write as \begin{equation}\label{d}\Fc \left(\tilde f(\theta,\rho)e^\rho\right)(\hat{\theta},\hat{\rho})=\frac{\Fc \left((\tilde{f}e^\rho)*B_3\right)(\hat{\theta},\hat{\rho})}{\Fc B_3(\hat{\theta},\hat{\rho})},\end{equation} where $B_3$ is the cubic (cardinal) B-spline. Here we only consider frequencies $(\hat{\theta},\hat{\rho})$ in a rectangle $L$, where $\left| \Fc B_3(\hat{\theta},\hat{\rho}) \right| $ does not become too small.

By using the coordinate transformations \eqref{f}, \eqref{transformT} and \eqref{transformP} let
\begin{equation}
\begin{pmatrix}\varphi(t,x)\\\eta(t,x)\end{pmatrix}=\mathsf{P_1T}\begin{pmatrix}t^2\\x^2\end{pmatrix}
\end{equation}
be log-polar coordinates that correspond to the coordinates $(t^2,x^2)$ in the time-offset domain. In these coordinates \eqref{d} takes the form
\begin{equation}\label{eq:B3tx}\begin{aligned}
\Fc \left(f(\theta,\rho)e^\rho\right)(\hat{\theta},\hat{\rho})=\frac{\Fc\left(\iint\! \frac{f(t,x)}{2x}e^\eta \mathbf{J}(\!t,x\!)B_3\!\left(\theta\!-\!\varphi,\rho\!-\!\eta\right) dtdx\right)(\hat{\theta},\hat{\rho})}{\Fc B_3(\hat{\theta},\hat{\rho})} \end{aligned}
\end{equation}
where the division by $2x$ is related to the transformation \eqref{f}. However, the Jacobian determinant $\mathbf{J}(t,x)=\left|\frac{\partial (\varphi,\eta)}{\partial (t,x)}\right|$ is a easily seen to consist of smooth bounded functions multiplied with $2x$ (coming from the derivative of $x^2$), which cancels out this seeming singularity at $x=0$. Subsequently the integrals and the Fourier transforms above can be well approximated using the trapezoidal rule and FFT for approximative evaluation of $\Fc$. If $(t_j,x_k)$ are regular sampling points in the time-offset domain, we introduce the approximation to \eqref{eq:B3tx} by 
\begin{equation}\label{eq:B3tx2}\begin{aligned}
g(\hat{\theta},\hat{\rho})=c\frac{\Fc\left(\sum_{j,k} \frac{f(t_j,x_k)}{2x_k}e^{\eta(t_j,x_k)} \mathbf{J}(\!t_j,x_k\!)B_3\!\left(\theta\!-\!\varphi(t_j,x_k),\rho\!-\!\eta(t_j,x_k)\right) \right)(\hat{\theta},\hat{\rho})}{\Fc B_3(\hat{\theta},\hat{\rho})} \end{aligned},
\end{equation}
where $c$ is a constant related to the sampling intervals. This approximation is then accurate for values of $(\hat{\theta},\hat{\rho})$ in the rectangle $L$ mentioned above.

As outlined in the previous section, this allows us to efficiently compute approximations of \eqref{finiteconv} on a regular lattice in the log-polar coordinate system $(\theta,\rho)$ via the formula \begin{equation}\label{lok}\Rad_\text{lp}\tilde{f}(\theta,\rho)\approx\cos(\theta)\Fc^{\text{-}1}\Big( g(\hat{\theta},\hat{\rho})\cdot\Fc \zeta(\hat{\theta},\hat{\rho})\Big)(\theta,\rho).\end{equation} The final interpolation from the log-polar $(\theta,\rho)$ lattice to the Radon $(\tau,q)$ lattice can be done by using cubic B-splines and a slight modification of \eqref{lok}. Here, we again employ the transformations (\ref{transformT}-\ref{transformP}). In this case, let
\begin{equation}
\begin{pmatrix}\varphi(\tau,q) \\ \eta(\tau,q)\end{pmatrix}=\mathsf{P}_2\mathsf{S}\begin{pmatrix}\tau^2 \\ q^2\end{pmatrix}
\end{equation}
be the log-polar coordinates that correspond to the coordinates $(\tau^2,q^2)$ in the Radon domain.
The interpolation from the log-polar $(\theta,\rho)$ lattice to the Radon $(\tau,q)$ lattice can then be done by using (a discrete version of)
\begin{align}\label{eq:B3tauq}
&\Radh f(\tau,q)=\cos(\varphi)\int_{-\frac{\beta}{2}}^{\frac{\beta}{2}}\int_{\log(a_r)}^{0} \Big(
\Fc^{-1}\left(\frac{\chi_L g(\hat{\theta},\hat{\rho}) \Fc{\zeta}(\hat{\theta},\hat{\rho}) }{\Fc B_3(\hat{\theta},\hat{\rho})}\right)(\theta,\rho)\Big) B_3\left(\varphi-\theta,\eta-\rho \right)  d\rho d\theta,
\end{align}
where $\chi_L$ denotes the characteristic function of the set $L$.

Numerical evaluation of the approximations \eqref{eq:B3tx} and \eqref{eq:B3tauq} appear to be well-suited for parallel computations, particularly on GPUs. For FFT we make use of the high-performance cuFFT library, efficient GPU kernels can be constructed for the smearing operations and for the vector multiplications. The discrete version of the operator $\Radh$, as explained in the previous sections, will be denoted by $R_\text{h}$.

\section{Reconstruction techniques}\label{rec_tec}
The adjoint operator for the hyperbolic Radon transform $R_\text{h}^*$ is defined by using the inner product equality
\begin{equation}\label{inner}
\langle R_\text{h} f,g\rangle=\langle f,R_\text{h}^* g\rangle,
\end{equation}
for arbitrary $f$ and $g$.
The operator is easy to construct by using the approach with switching to log-polar coordinates, essentially by reversing the order of the operations. With the adjoint operations at hand, one can consider iterative methods for representing $f$ by sparse sums of hyperbolic wave events, and related interpolation and reconstruction techniques. A popular such method is based on the soft thresholding method for obtaining sparse representations proposed in \cite{daubechies2004iterative}. In this setting it means to consider the minimization of
\begin{equation} \label{R_ell1}
\|R_\text{h}^*g-f\|_2^2+\mu\|g\|_1,
\end{equation}
for some choice of sparsity parameter $\mu$. 

By a simple modification of Theorem 3.1 in \cite{daubechies2004iterative}), this minimization problem is solved by the iterations 
\begin{equation}\label{iters}
g^n=\mathbf{S}_{c^2\mu}(g^{n-1}+c^2R_\text{h}(f-R_\text{h}^* g^{n-1})),\quad n=1,2,\dots,
\end{equation}

where $g^0$ is arbitrary, $c$ is a positive constant such that $c\|R_\text{h}\|<1$, and $\mathbf{S}_\mu$ is a soft-thresholding function defined as
\begin{equation}
\mathbf{S}_\mu(v)=\begin{cases}
v+\frac{\mu}{2}, & \mbox{if } v\le -\frac{\mu}{2}, \\
0, & \mbox{if } |v|<\frac{\mu}{2},\\
v-\frac{\mu}{2}, & \mbox{if } v\ge \frac{\mu}{2}.
\end{cases}
\end{equation}

To perform interpolation in the case of missing data, let $S$ be a subset of the $(t_j,x_k)$ grid where we do have measurements of $f$. We are then interested in minimizing \begin{equation}
\sum_{(t_j,x_k)\in S}(R_\text{h}^*g-f)^2(t_j,x_k)+\mu\|g\|_1,
\end{equation}
which, defining $f$ to be 0 where data is missing, is solved by the iteration 
\begin{equation}\label{itersmd}
g^n=\mathbf{S}_{c^2\mu}(g^{n-1}+c^2R_\text{h}(f-\chi_S R_\text{h}^* g^{n-1})),\quad n=1,2,\dots.
\end{equation}
Here $\chi_S$ is the characteristic function of $S$. Again, this scheme is efficiently evaluated using the fast implementation of 
$R_\text{h}$ explained in the previous section. 

\section{Discretization}\label{ch_disc}
In this section we derive guidelines for how to choose discretization parameters. For simplicity, we assume to work with regular sampling in the time-offset, and in the Radon domain; but the log-polar-based method can be easily generalized for unequally-spaced grids in these two domains.

In order to apply FFTs, samples in log-polar coordinates $(\theta,\rho)\in \left[-\betah,\betah \right]\times[\log a_r, 0]$ must be chosen on an equally spaced grid. By using coordinate conversions for the log-polar setup given by $$\begin{pmatrix}\varphi(t,x)\\\eta(t,x)\end{pmatrix}=\mathsf{P_1T}\begin{pmatrix}t^2\\x^2\end{pmatrix}.$$ In order to maintain accurate interpolation, we choose the sample spacing in $\theta$ and $\rho$ with respect to the largest distance between sample points in the $\varphi$ and $\eta$ variables, i.e.,
\begin{equation}\begin{aligned}\label{lpsamples}
&\Delta \theta \ge \max_{t_j,x_k}\left(\left|\varphi(t_j,x_k)-\varphi(t_j+\Delta t,x_k)\right|,\left|\varphi(t_j,x_k)-\varphi(t_j,x_k+\Delta x)\right|\right),\\
&\Delta \rho \ge \max_{t_j,x_k}\left(\left|\eta(t_j,x_k)-\eta(t_j+\Delta t,x_k)\right|,\left|\eta(t_j,x_k)-\eta(t_j,x_k+\Delta x)\right|\right).\end{aligned}
\end{equation}
This choice will determine the log-polar frequency range that can be covered, which in turn determines the resolution in the  $(\tau,q)$ (Radon) domain. 
The quadratic behavior in the time sampling can be fairly well described in terms of the log-polar sampling, as long as time range is not too large. In the case of large time ranges, it can be beneficial to split the split the time-offset and Radon domains in parts and consider the log-polar Radon transform for each of these parts, in order to avoid too large differences in sample densities.

\begin{figure}
	\centering 	
	\begin{subfigure}[b]{0.32\textwidth} \includegraphics[width=1\textwidth]{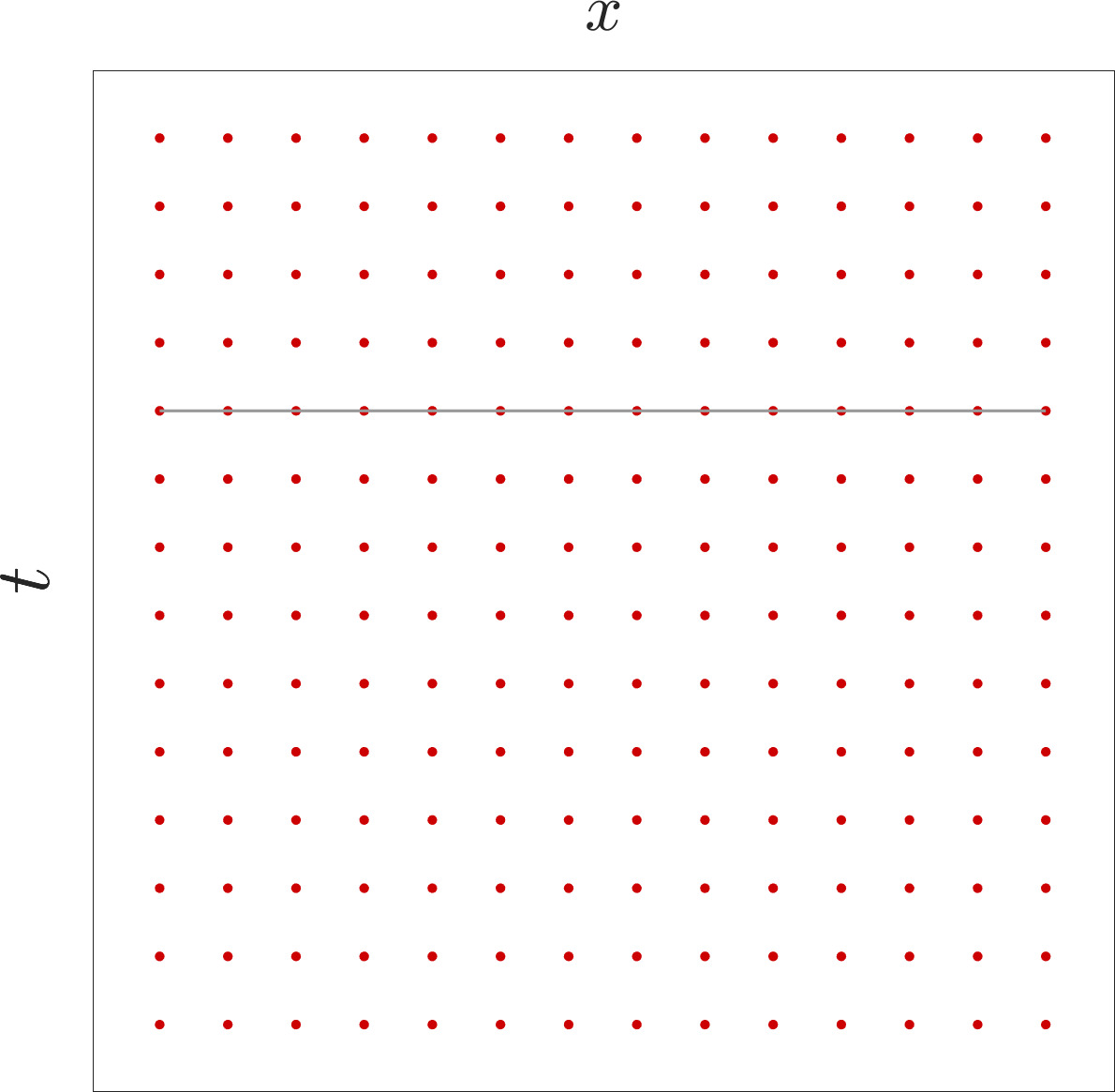}
		\caption{Grid $(t,x)$ with splitting in $t$ coordinate}
	\end{subfigure}\hspace{5pt}
	\begin{subfigure}[b]{0.32\textwidth} \includegraphics[width=1\textwidth]{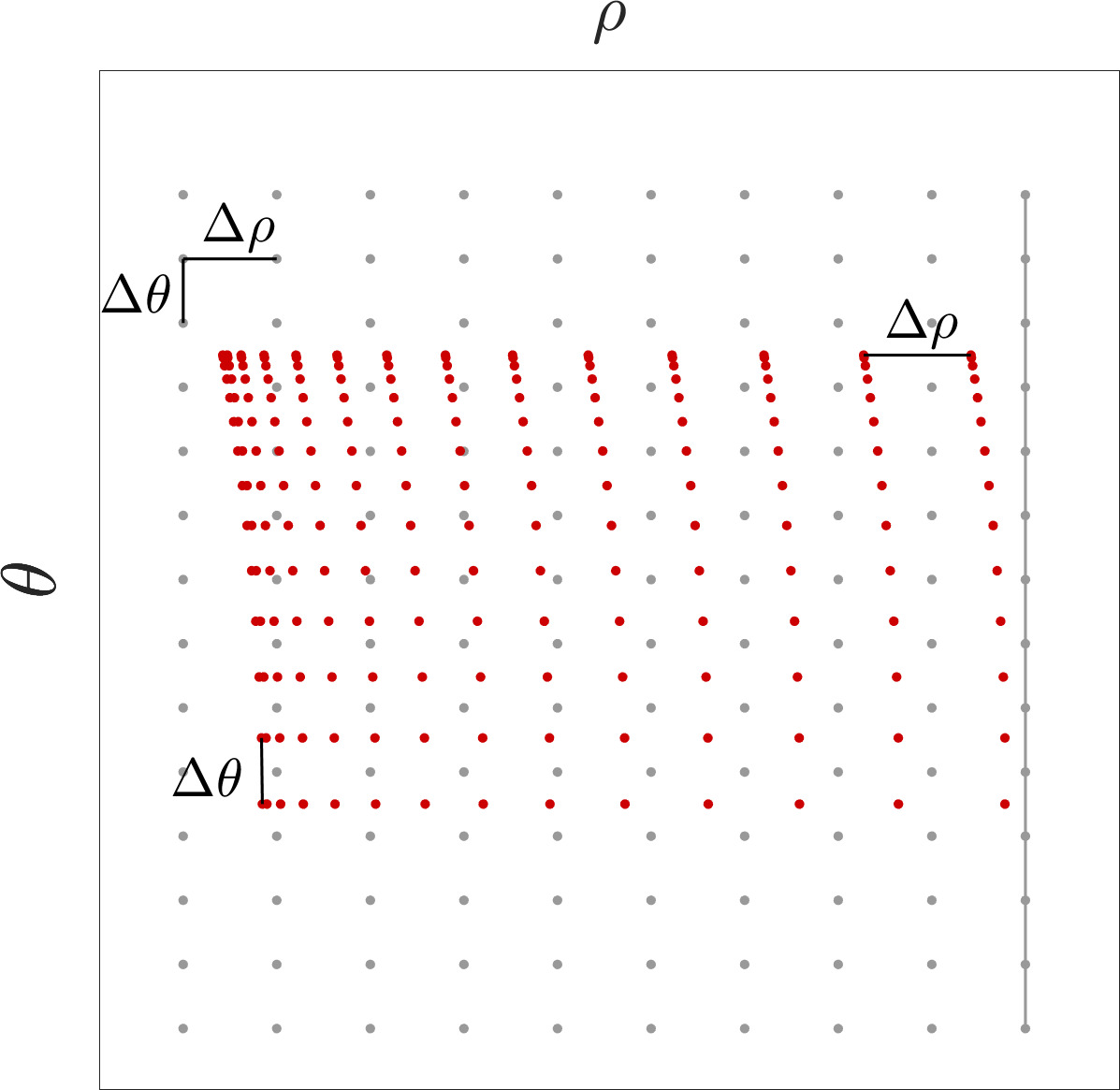}
		\caption{Grid $(t,x)$ in log-polar coordinates}	
	\end{subfigure}\hspace{5pt}
	\begin{subfigure}[b]{0.32\textwidth} \includegraphics[width=1\textwidth]{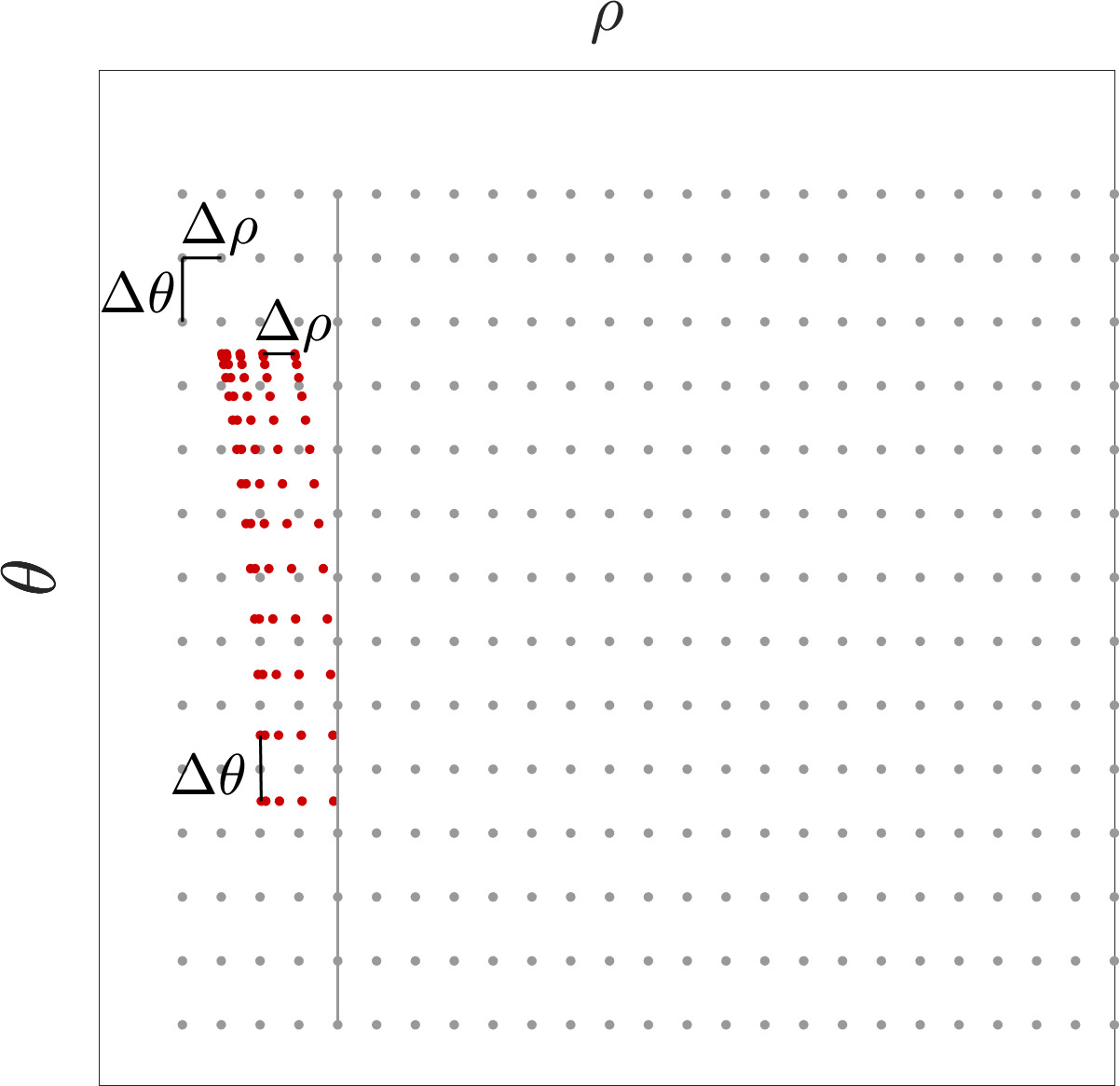}
		\caption{Part in log-polar coordinates for small values of $t$ }
	\end{subfigure}
	\caption{Grids for conversion between time-offset and log-polar coordinates. Splitting in $t$ variable.}
	\label{fig:lpsamples}	
\end{figure}

For instance, for small values of $\tau$ the grid for the Radon domain becomes more dense and samples $\Delta \theta,\Delta \rho$ should be chosen to be smaller. Suppose that we have already rescaled $f$ according to \eqref{hypnorm}, and note that the function is 0 until the arrival of the first event at $\tau_{min}$. We may then split the integral in the following way
\begin{equation}\label{tsplit}\begin{aligned}
&\Radh f(\tau,q)=\int_{\tau_{min}}^{1}\int_{0}^{1} f(t,x)\delta(t-\sqrt{\tau^2+q^2 x^2})dx dt=\\&\int_{\tau_{min}}^{a}\int_{0}^{1} f(t,x)\delta(t-\sqrt{\tau^2+q^2 x^2})dx dt + \int_{a}^{1}\int_{0}^{1} f(t,x)\delta(t-\sqrt{\tau^2+q^2 x^2})dx dt
\end{aligned}
\end{equation}
for some $a$ between $\tau_{min}$ and $1$. For numerical evaluation of the first integral by using the log-polar based method samples in $\theta,\rho$ determined according to \eqref{lpsamples} become more dense, see Figure \ref{fig:lpsamples} for a schematic description. The red dots in Figure \ref{fig:lpsamples}b indicate log-polar samples after conversion to discrete coordinates in the $(t,x)$ domain illustrated in Figure \ref{fig:lpsamples}a. Equally spaced samples in the log-polar domain (gray dots) are chosen with respect to maximal distances between points \eqref{lpsamples}. Figure \ref{fig:lpsamples}c demonstrates samples in the log-polar domain corresponding to small values of $t$ (located above the gray line in Figure \ref{fig:lpsamples}a). The splitting procedure is not computationally intensive and can be applied several times to achieve accuracy for small values of $\tau$.

\section{Accuracy and performance tests}
\begin{figure}[t!]
	\centering	
	\begin{subfigure}[b]{0.425\textwidth} \includegraphics[width=1\textwidth]{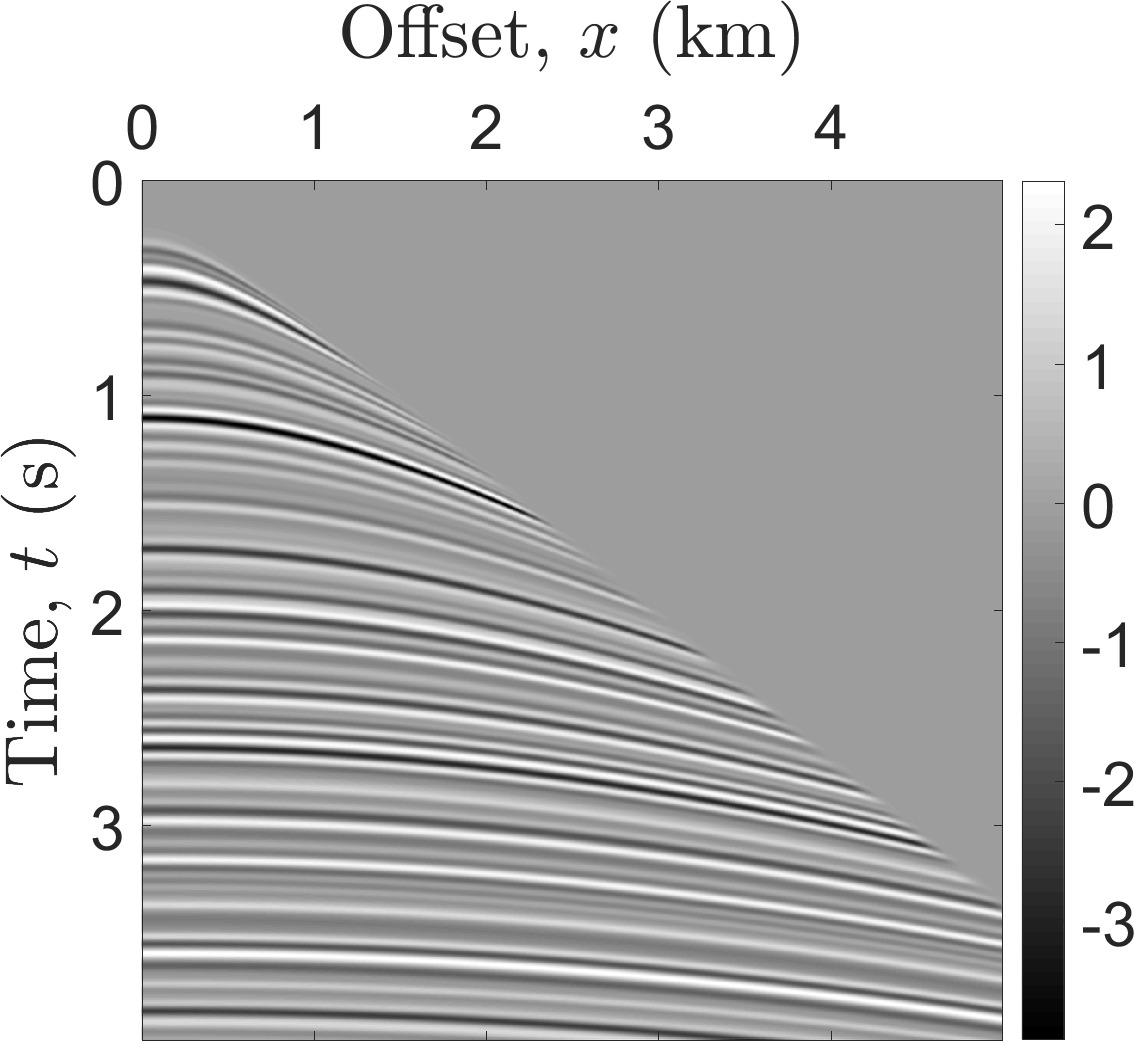}
		\caption{CMP gather}
	\end{subfigure}	\hspace{5pt}
	\begin{subfigure}[b]{0.425\textwidth} \includegraphics[width=1\textwidth]{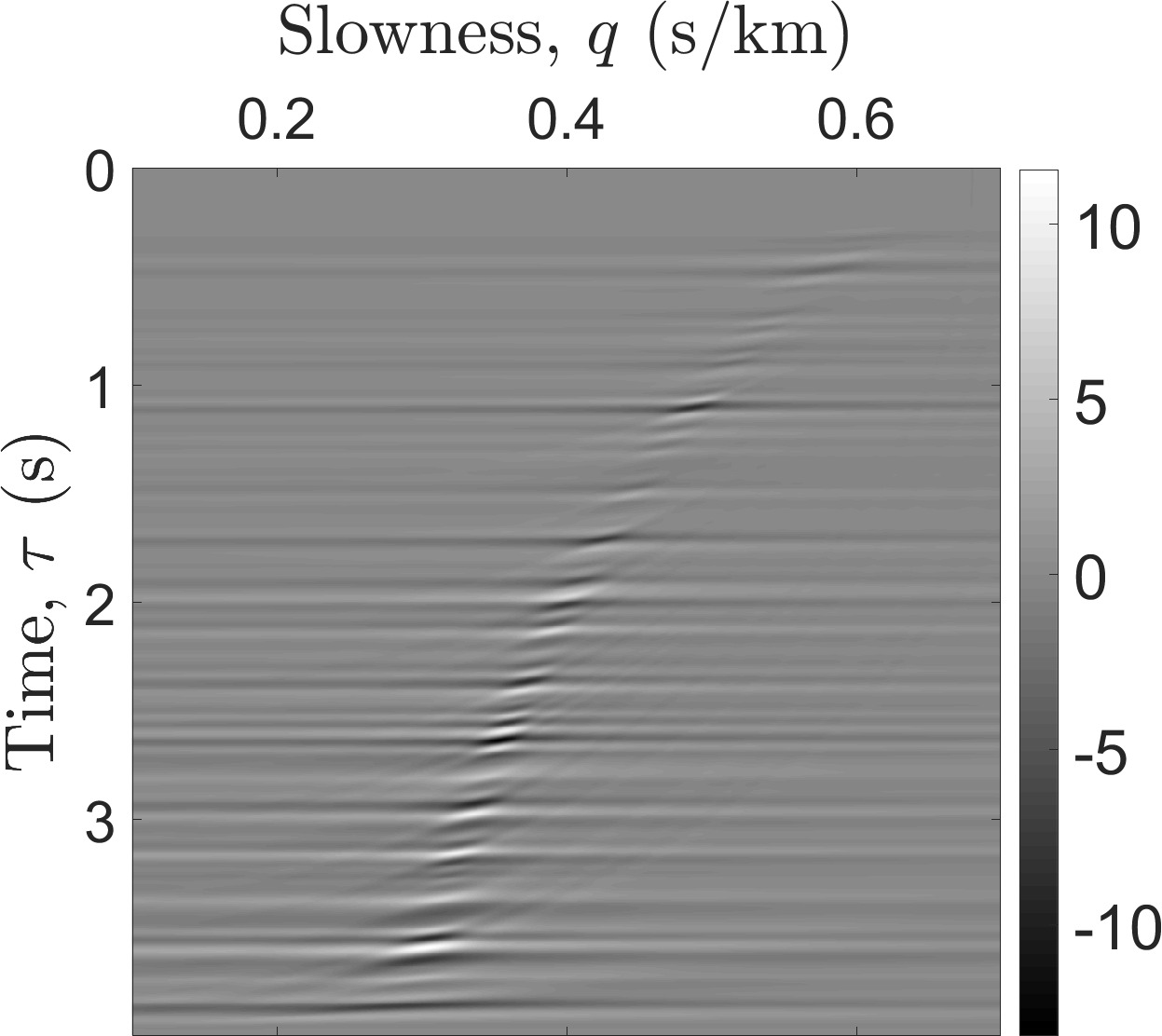}
		\caption{Hyperbolic Radon transform}
	\end{subfigure}		
	
	\vspace{5pt}
	
	\begin{subfigure}[b]{0.425\textwidth} \includegraphics[width=1\textwidth]{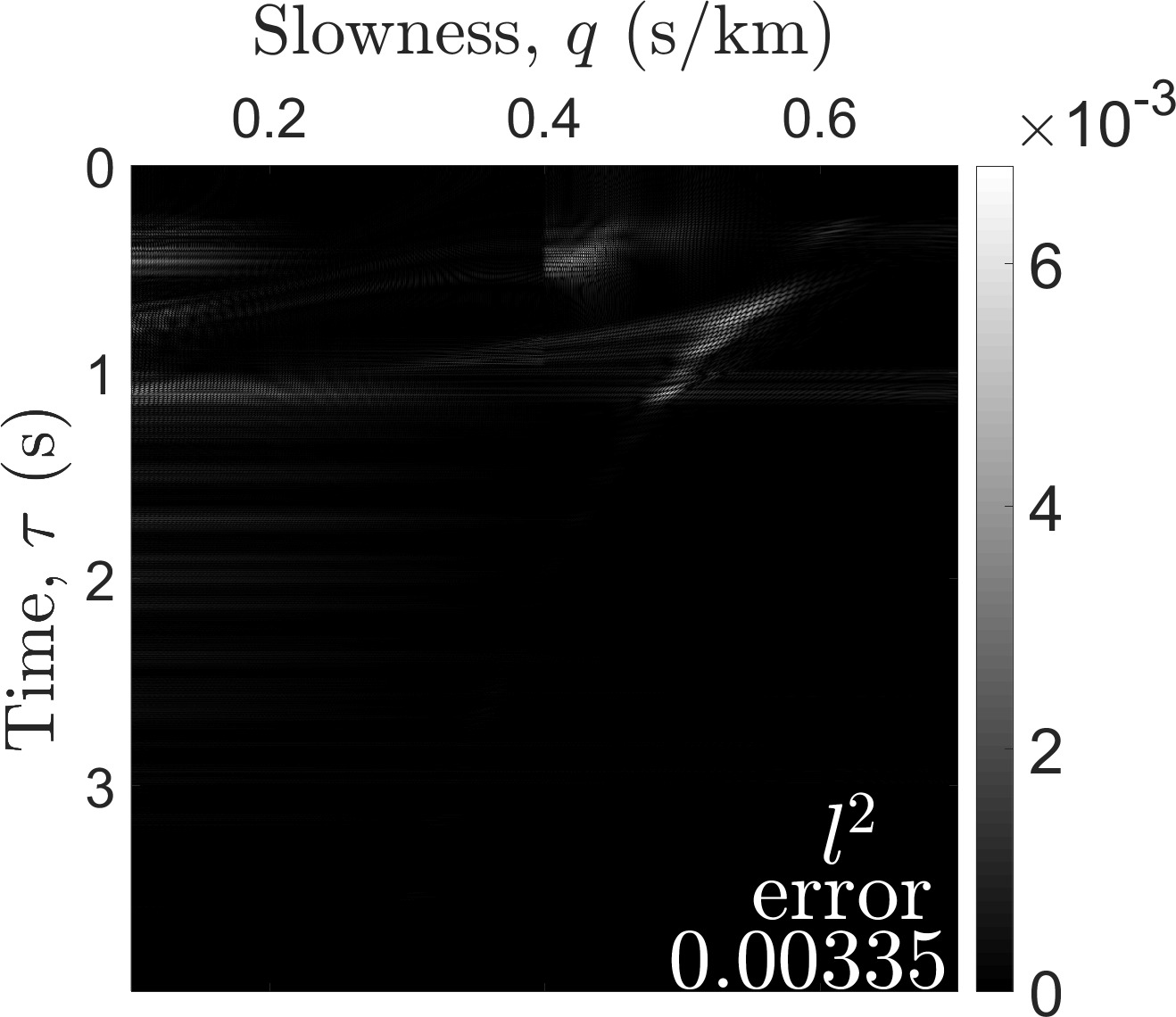}
		\caption{Normalized error, Log-polar-based method}
	\end{subfigure}	\hspace{5pt}
	\begin{subfigure}[b]{0.425\textwidth} \includegraphics[width=1\textwidth]{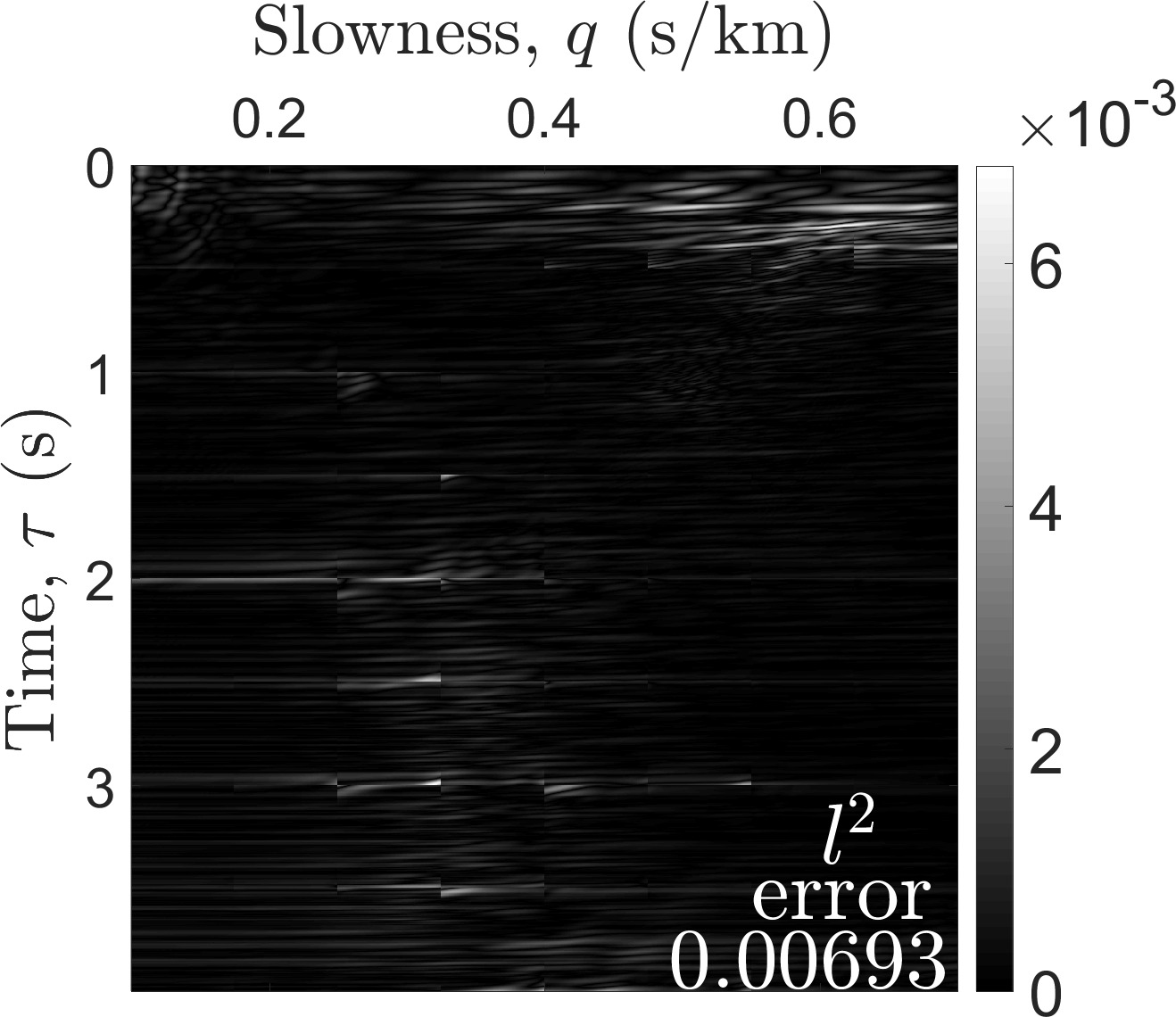}
		\caption{Normalized error,
			Fast butterfly algorithm ($N=64$, $q_i=9$)}
	\end{subfigure}	
	\caption{Hyperbolic Radon transform. Corresponding normalized errors compared to direct summation over hyperbolas. }
	\label{fig:cmpsyn}
\end{figure}

For the sake of quality comparisons, we perform the same tests as the ones presented by \cite{hu2013fast} for the fast hyperbolic Radon transform based on fast butterfly algorithms. The method of \cite{hu2013fast} is available in the open source software package Madagascar \cite{fomel2013madagascar}. The synthetic CMP gather (Figure \ref{fig:cmpsyn}a) was used as a reference for making comparisons. As a reference method, we use a standard C implementation of the direct summation given by \eqref{Rhdef2}. Here cubic interpolation is used for the interpolation in time.

The fast butterfly algorithm has several parameters for controlling efficiency and accuracy, for details we refer to the pages 5, 6 in \cite{hu2013fast}. The parameter $M$ ($N$ in the paper) is of the order of the maximum value of the phase function $|\Phi(\mathbf{x},\mathbf{k})|$ used for the approximation; and parameters $q_{k_1},q_{k_2},q_{x_1},q_{x_2}$ control the number of Chebyshev points. According to the results presented in \cite{hu2013fast}, the set of parameters $(q_i=9, M=64)$ shows an accuracy level of about $O(10^{-3})$ for images of size $1000\times1000$. We performed tests for $N_t=N_x=N_\tau=N_q=N$ where $N$ was chosen as different powers of 2, and for obtaining an approximate accuracy level of $\Oc(10^{-3})$ we used $(q_i=9, M=N/16)$ in accordance with the tests conducted in \cite{hu2013fast}.

Normalized errors compared to direct summation over hyperbolas for the log-polar-based and for the fast butterfly algorithm are demonstrated in Figure \ref{fig:cmpsyn}c and Figure \ref{fig:cmpsyn}d, respectively. The figures show that the two methods have the same order of errors. The errors for the log-polar-based method are mostly observed in the region of small time intercept ($\tau$) and high values of slowness ($q$). These accuracy problems can be reduced by additional splittings of the integral for the hyperbolic Radon transform, similar to the one suggested in the expression \eqref{tsplit}. To be concrete, the presented results for the log-polar-based method were obtained after one splitting in the time variable, and one splitting in the slowness variable.

Table \ref{tab:times} demonstrates the computational times for the fast butterfly algorithm; for the CPU and the GPU versions of the log-polar-based method, respectively; and for the direct summation over hyperbolas. The table confirms the complexity of the proposed method and shows that a substantial performance gain is obtained by using GPUs. It is common in GPU computing that time to copy data between host and device memory constitute an essential part of the total computational costs (for our tests it takes $\approx 30\%$ of total time). This time can be neglected in the case of using iterative schemes since it is then possible to keep all data in the GPU memory. For the tests performed, we used a standard desktop with an Intel Core i7-3820 processor and NVIDIA GeForce GTX 970 video card with PCI Express x16 graphic interface. All computations were performed in single precision. We note that parallel versions of butterfly algorithms have been described in \cite{poulson2014parallel}, but we use the single core implementation described in \cite{hu2013fast} to make sure that the computational times are in accordance with the results reported in  \cite{hu2013fast}.

\begin{figure}[t!]	
	\centering 	
	\includegraphics[width=0.42\textwidth]{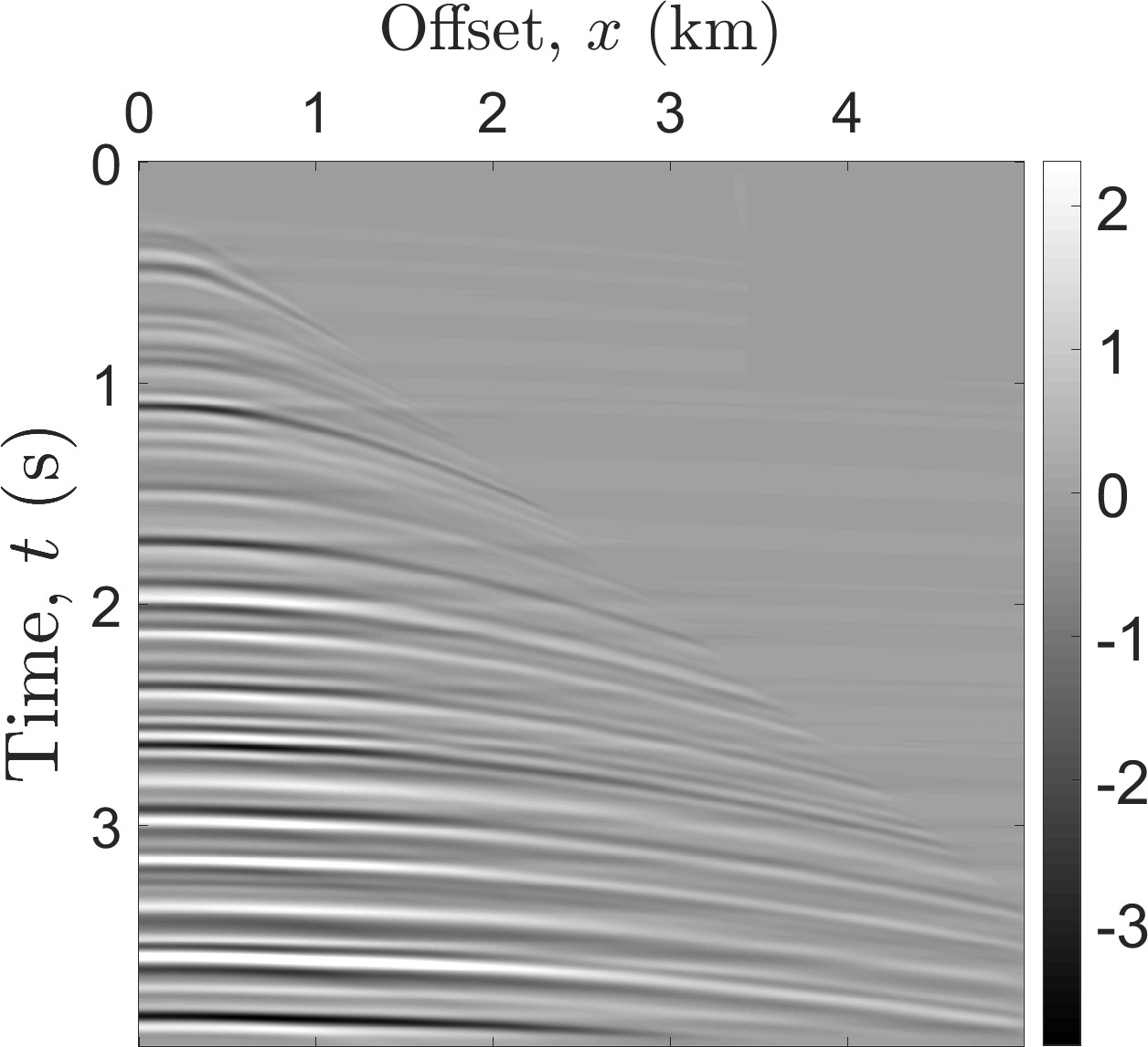}\hspace{5pt}
	\includegraphics[width=0.42\textwidth]{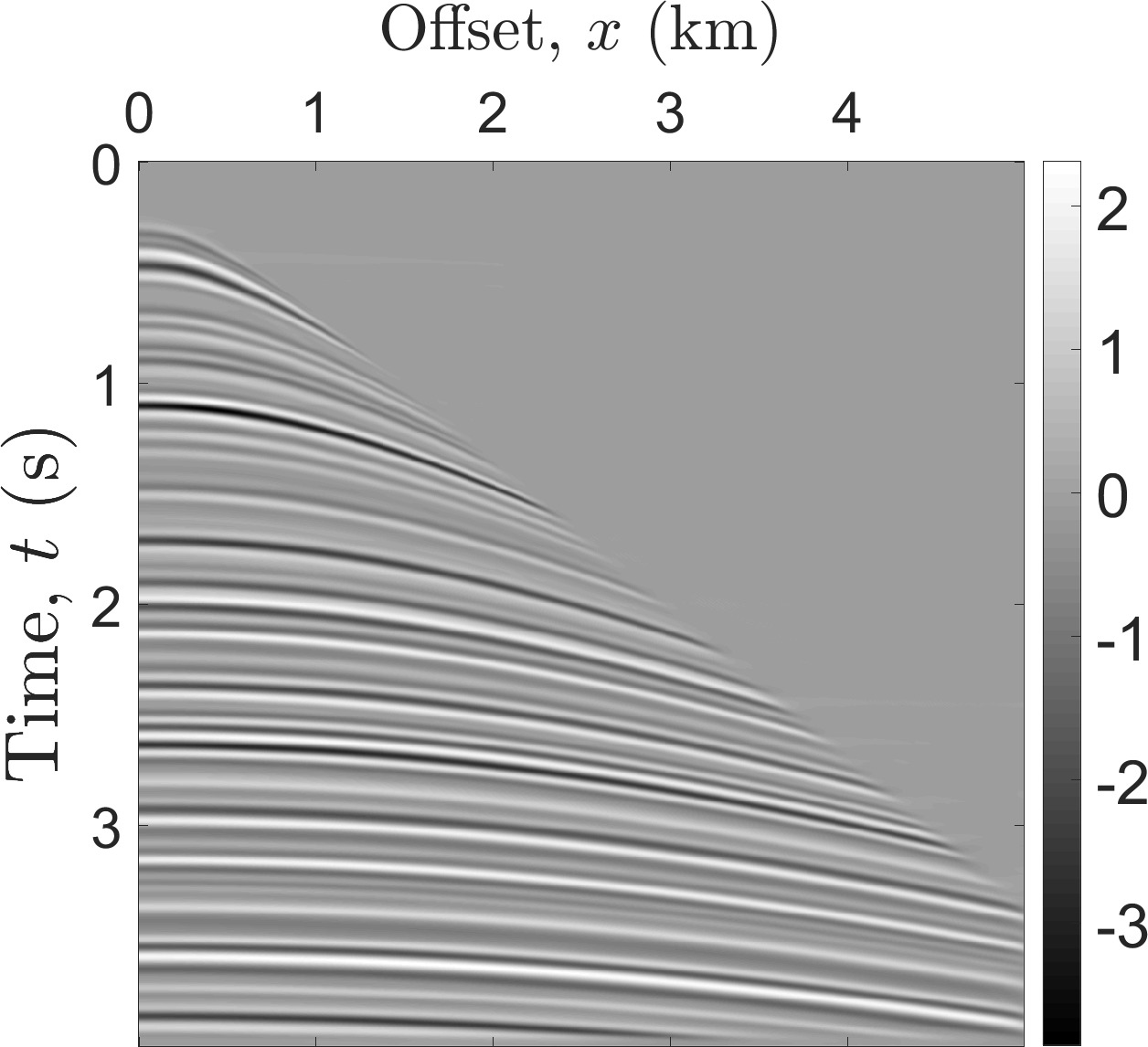}
	\caption{Output of the adjoint operator for the hyperbolic RT (left), and the result from using 30 soft-thresholding iterations from \protect\eqref{iters} (right).}	
	\label{fig:adj}	
\end{figure}

In Figure \ref{fig:adj} we show the output of the adjoint operator for the hyperbolic Radon transform, as well as $\ell^1$ regularized reconstruction given by \eqref{R_ell1}. The proposed algorithm passes the inner product test \eqref{inner} with a relative error $\Oc(10^{-5})$. The iterative reconstruction demonstrates good quality (compare  figures \ref{fig:cmpsyn}a and \ref{fig:adj}b).

\begin{table}[h!]\centering\label{tab:timefwd}
	\caption{Computational time (in sec) for the hyperbolic Radon transform via direct summation over hyperbolas, fast butterfly algorithm and via the log-polar-based method (CPU and GPU), speed-up compared to the direct summation.}\label{tab:times}
	\begin{tabular}{|c|c|c|c|c|c|c|c| } 		
		\hline
		& \thead{Direct sums \\ CPU, 1 core} & \multicolumn{2}{c|}{\thead{Fast butterfly \\ CPU, 1 core}}& \multicolumn{2}{c|}{\thead{Log-polar \\ CPU, 8 cores}} & \multicolumn{2}{c|}{\thead{Log-polar \\ GPU}}  \\
		\hline
		N & time & time & speed-up & time & speed-up & time & speed-up\\
		\hline
		$2^{9}$ & 4.8e\Plus00 & 1.1e\Plus00 &4.3 & 3.3e\Minus02 &145.5 &2.6e\Minus03 &1828.1\\
		\hline
		$2^{10}$ & 4.0e\Plus01 & 4.5e\Plus00 &9.0 & 1.2e\Minus01 &344.4 &9.6e\Minus03 &4220.2\\
		\hline
		$2^{11}$ & 3.2e\Plus02 & 1.8e\Plus01 &17.8 & 4.7e\Minus01 &682.4 &3.5e\Minus02 &9018.3\\
		\hline
		$2^{12}$ & 2.5e\Plus03 & 7.3e\Plus01 &33.7 & 2.0e\Plus00 &1257.7 &1.4e\Minus01 &17036.3\\
		\hline
	\end{tabular}
\end{table}

\begin{table}[h!]\centering
	\caption{Computational time (in sec) for 64 soft-thresholding iterations.}\label{tab:timesiters}
	\begin{tabular}{|c|c|c|c| } 		
		\hline
		\multirow{ 2}{*}{N}  & \multirow{ 2}{*}{Total time} & \multicolumn{2}{c|}{Average time per iteration} \\
		\cline{3-4}
		&   & Forward operator & Adjoint operator  \\
		\hline
		$2^{9}$ & 8.1e\Plus00 & 1.9e\Minus03 & 2.3e\Minus03\\
		\hline
		$2^{10}$ & 3.1e\Plus01 & 7.4e\Minus03 & 8.7e\Minus03\\
		\hline
		$2^{11}$ & 1.2e\Plus02 & 2.9e\Minus02 & 3.2e\Minus02\\
		\hline
		$2^{12}$ & 4.6e\Plus02 & 1.1e\Minus01 & 1.3e\Minus01\\
		\hline
	\end{tabular}
\end{table}
Table \ref{tab:timesiters} shows computational times  using a GPU implementation of the proposed log-polar-based hyperbolic Radon transform, and 64 iterations of the iterative scheme \eqref{iters}. The table also contains times for single iteration of the forward and adjoint operators. One can see that in comparison to the GPU results in Table \ref{tab:times}, the times for the forward operator are lower due to limited number of host-device data transfers. For this scheme data was copied only for an initial guess $g^0$; the measured data $f$; and the final result.

\section{Applications}
In this section we mention some applications of the fast hyperbolic Radon transform. These are fairly standard, but the examples could be of practical interest due to the substantial computational speedup of the proposed implementation of the hyperbolic Radon transform. 

\subsection{Multiple attenuation.}
\begin{figure}[t!]
	\centering 	
	\begin{subfigure}[b]{0.425\textwidth} \includegraphics[width=1\textwidth]{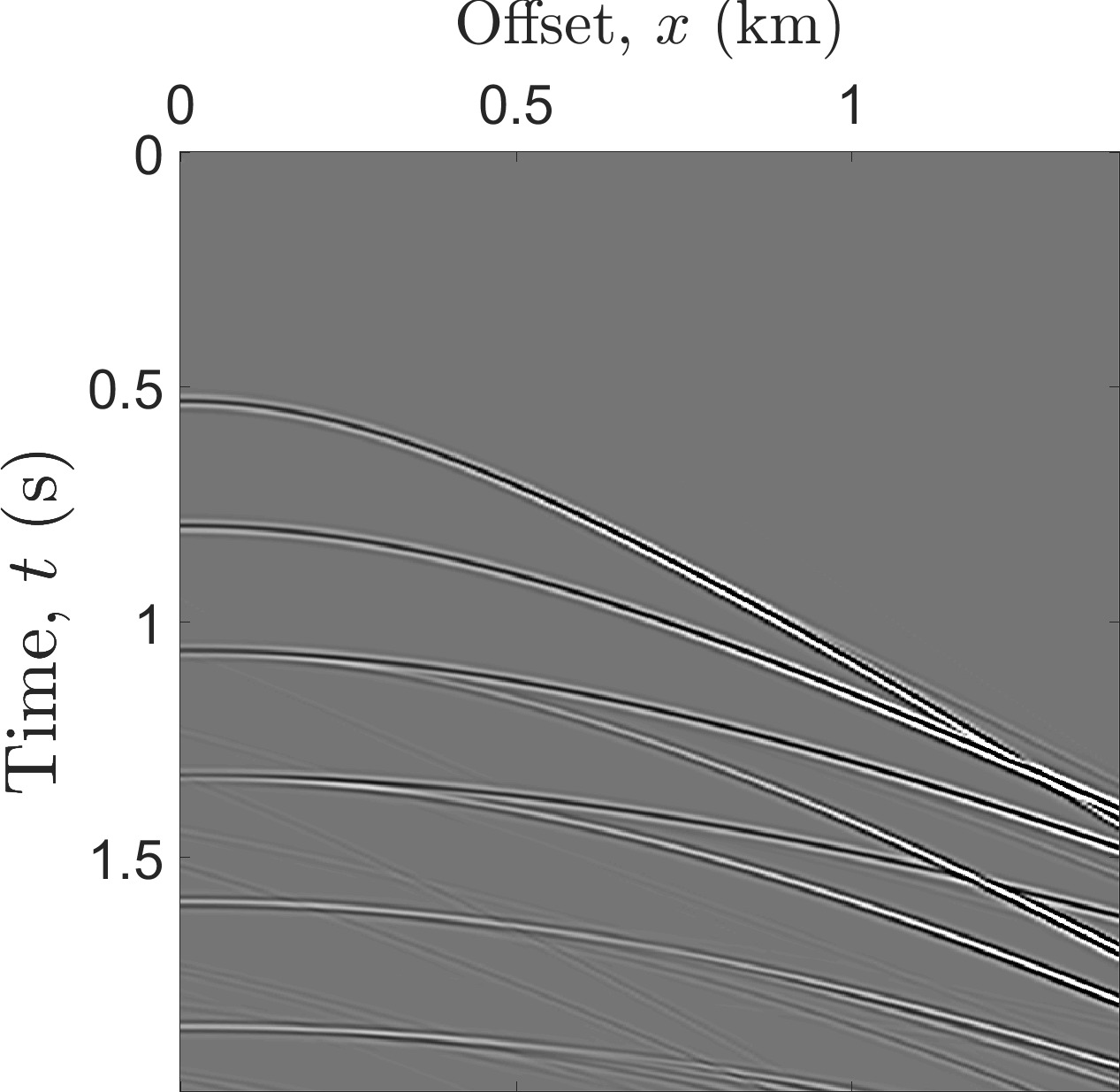}
		\caption{CMP gather}
	\end{subfigure}	\hspace{5pt}
	\begin{subfigure}[b]{0.425\textwidth} \includegraphics[width=1\textwidth]{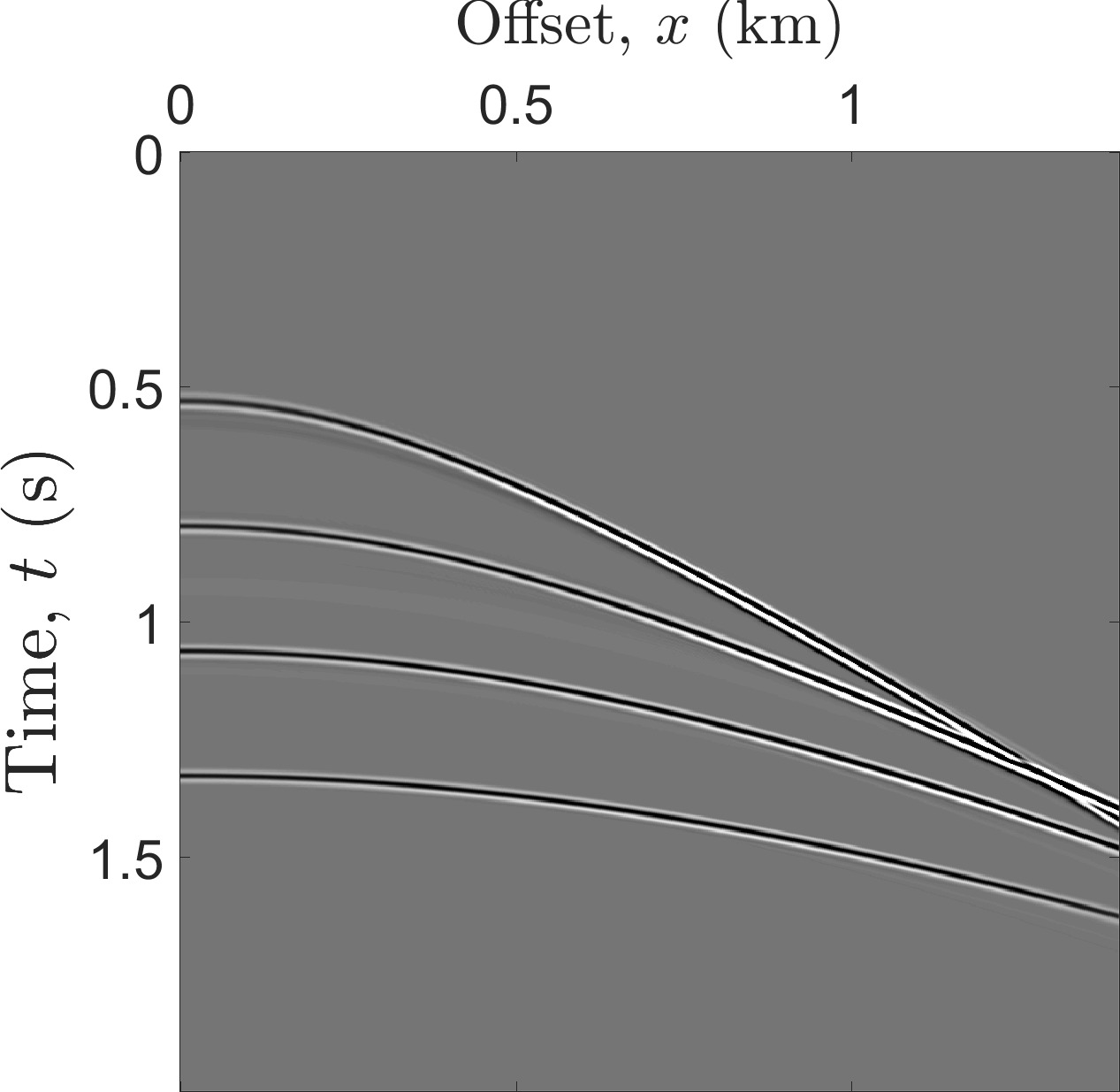}
		\caption{Reconstructed primaries}
	\end{subfigure}	
	
	\vspace{5pt}
	
	\begin{subfigure}[b]{0.425\textwidth} \includegraphics[width=1\textwidth]{cmp0iterstauqre.jpg}
		\caption{Hyperbolic Radon transform}
	\end{subfigure}	\hspace{5pt}
	\begin{subfigure}[b]{0.425\textwidth} \includegraphics[width=1\textwidth]{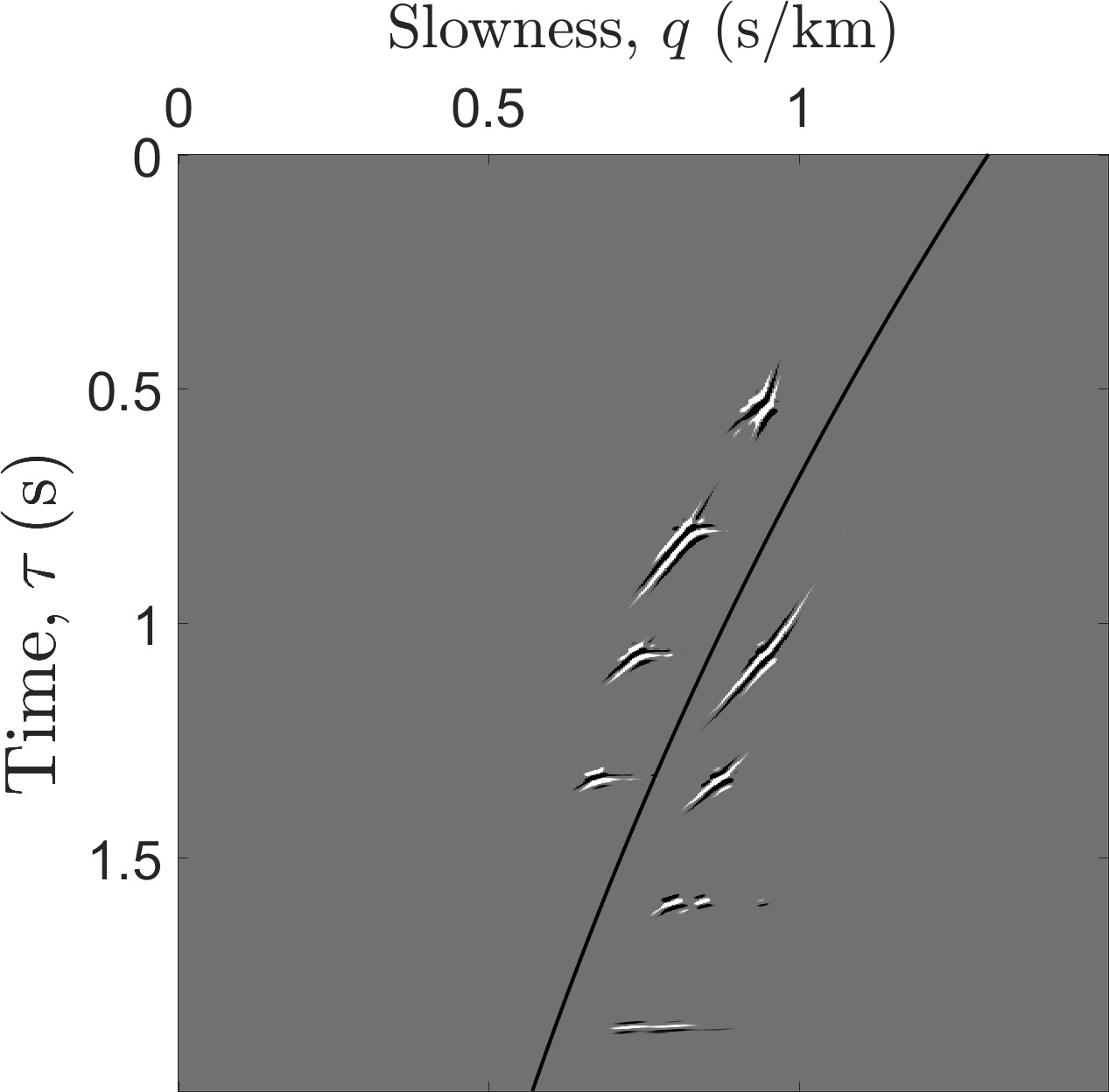}
		\caption{Sparse representation}
	\end{subfigure}	
	\caption{Multiple attenuation with 30 soft-thresholding iterations.}
	\label{fig:demulhyb}
\end{figure}
A well-known method for the attenuation of multiple reflections in CMP gathers is  based on conducting the attenuation in a Radon domain. Here, multiples and primaries can be separated due to their differences in moveout. We have tested method described in Chapter \ref{rec_tec} for the synthetic CMP gather in Figure \ref{fig:demulhyb}a.  Figure \ref{fig:demulhyb}c  illustrates the Radon data, and note that the primaries and multiples are difficult to separate. The corresponding result after using 30 soft-thresholding iterations is illustrated in Figure \ref{fig:demulhyb}d. The black line indicates the border between primaries and multiples, and Figure \ref{fig:demulhyb}b shows the reconstructed primaries.

\subsection{Interpolation.}
\begin{figure}[t!]
	\centering 	
	\begin{subfigure}[b]{0.425\textwidth} \includegraphics[width=1\textwidth]{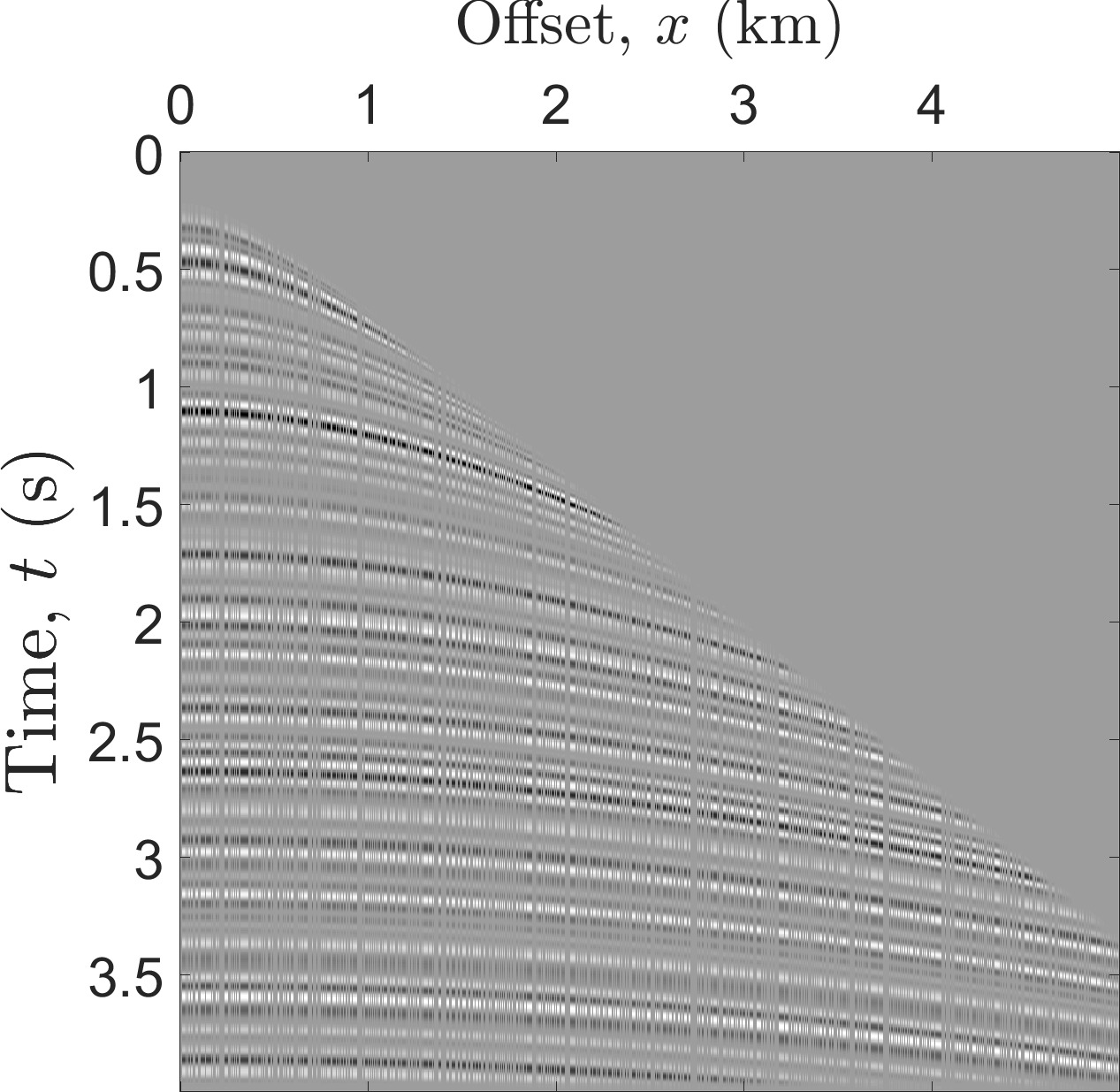}
		\caption{50\% missing traces}
	\end{subfigure}	\hspace{5pt}
	\begin{subfigure}[b]{0.425\textwidth} \includegraphics[width=1\textwidth]{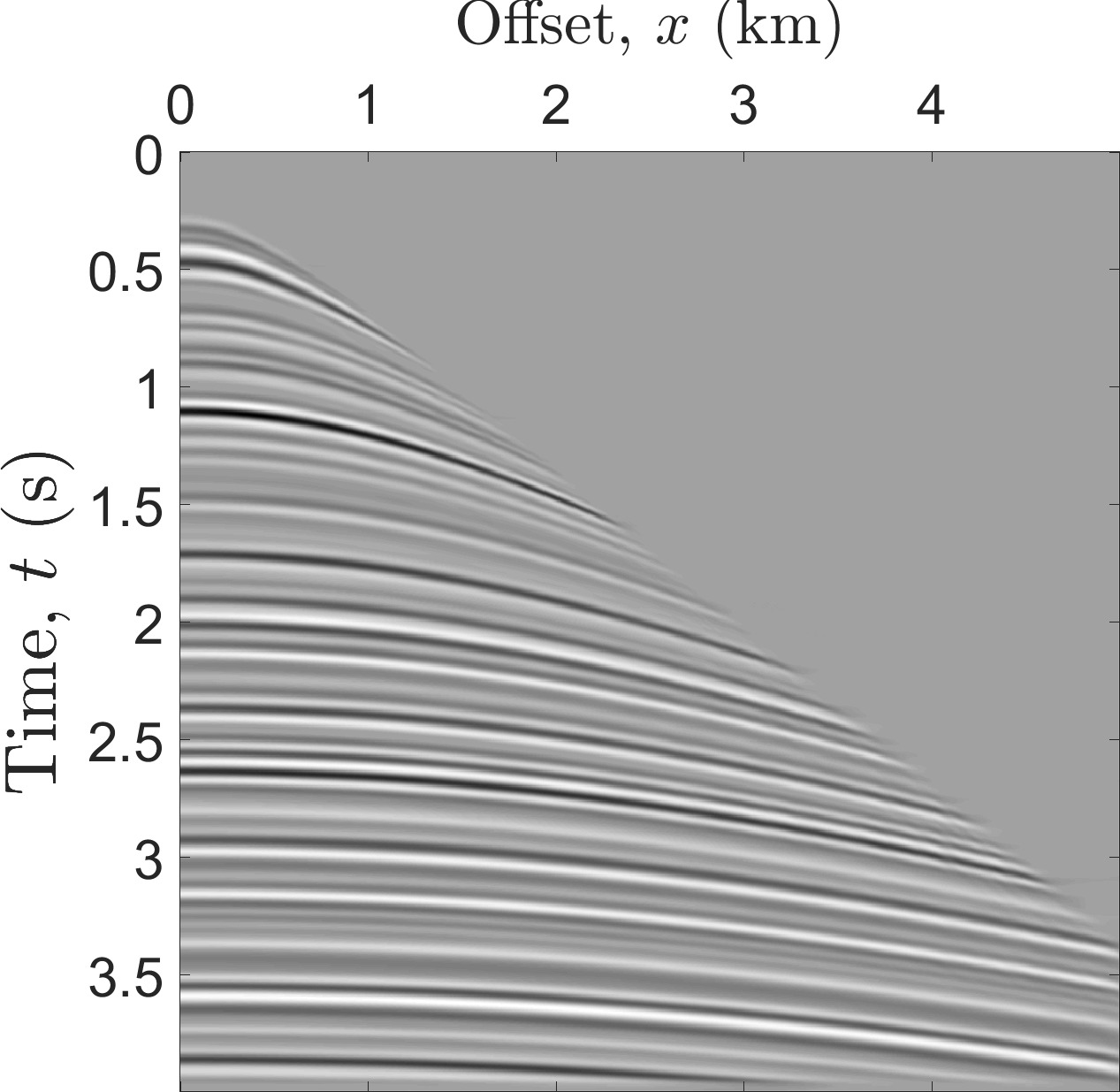}
		\caption{Reconstruction}		
	\end{subfigure}	
	\caption{Interpolation into missing traces with the soft-thresholding algorithm.}	
	\label{fig:cmpmissed}
\end{figure}
\begin{figure}[t!]
	\centering
	\begin{subfigure}[b]{0.425\textwidth}	
		\includegraphics[width=1\textwidth]{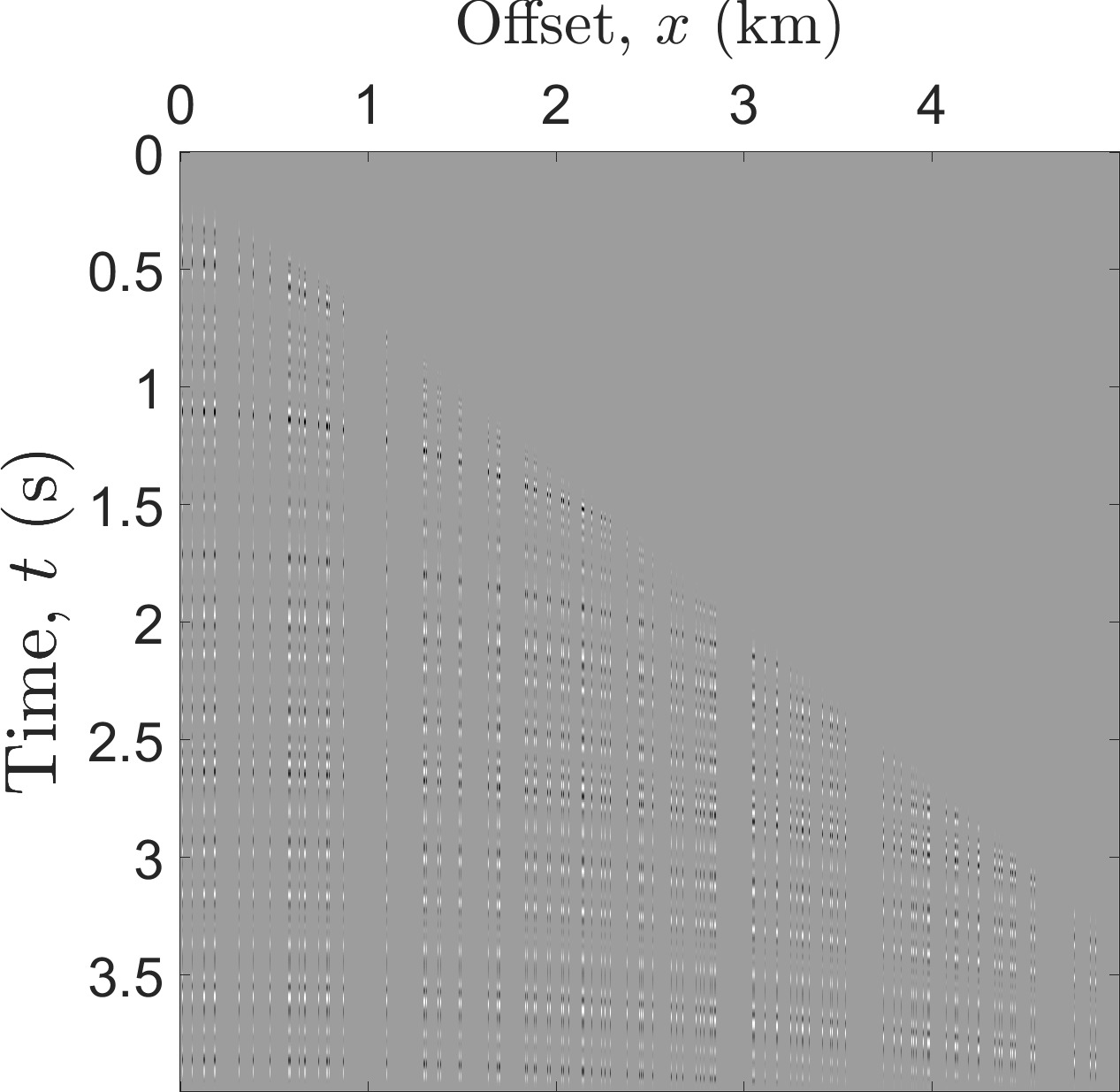}
		\caption{90\% missing traces}
	\end{subfigure}	
	\begin{subfigure}[b]{0.425\textwidth} \includegraphics[width=1\textwidth]{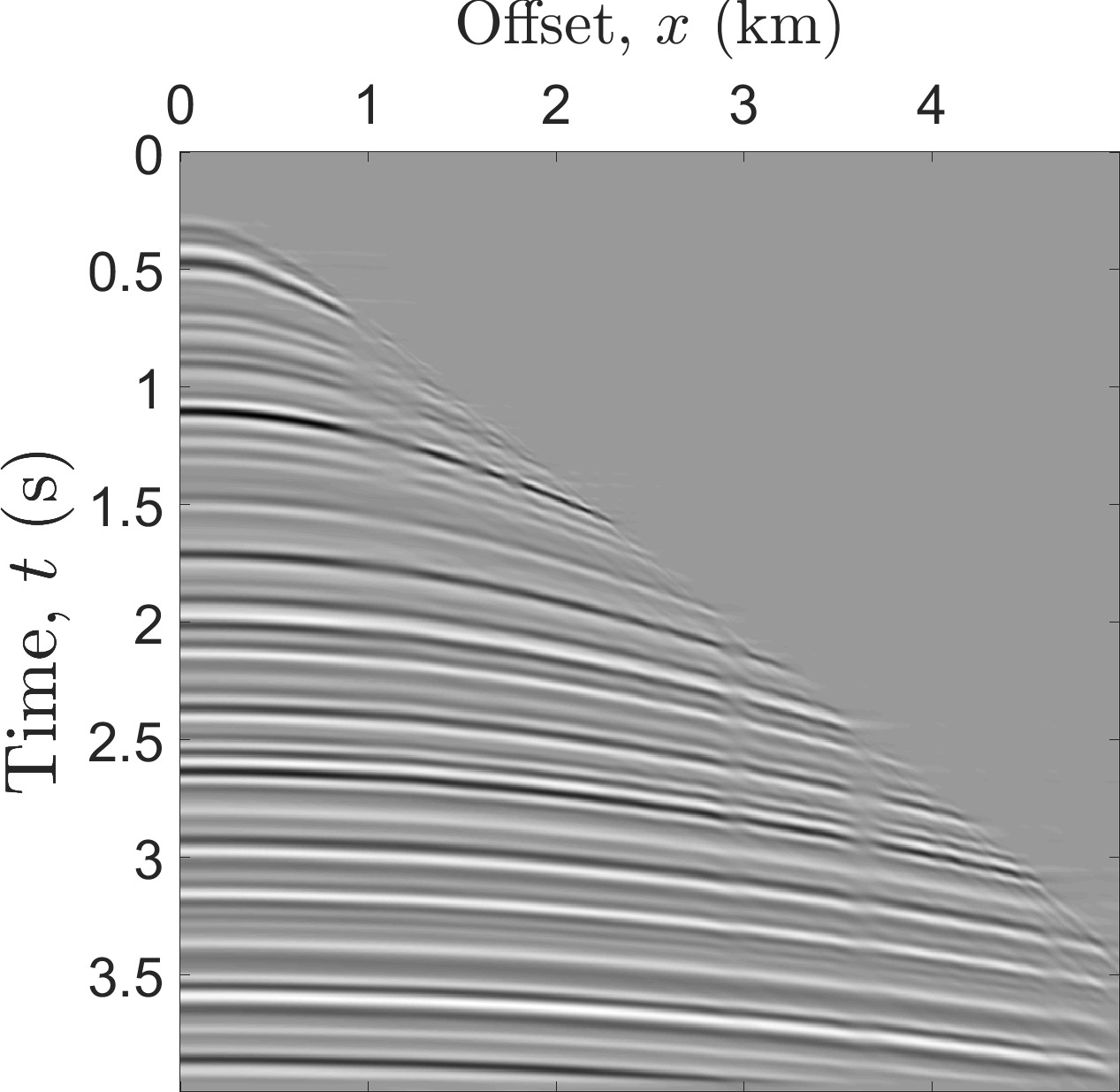}
		\caption{Reconstruction}		
	\end{subfigure}		
	
	\begin{subfigure}[b]{0.425\textwidth}	
		\includegraphics[width=1\textwidth]{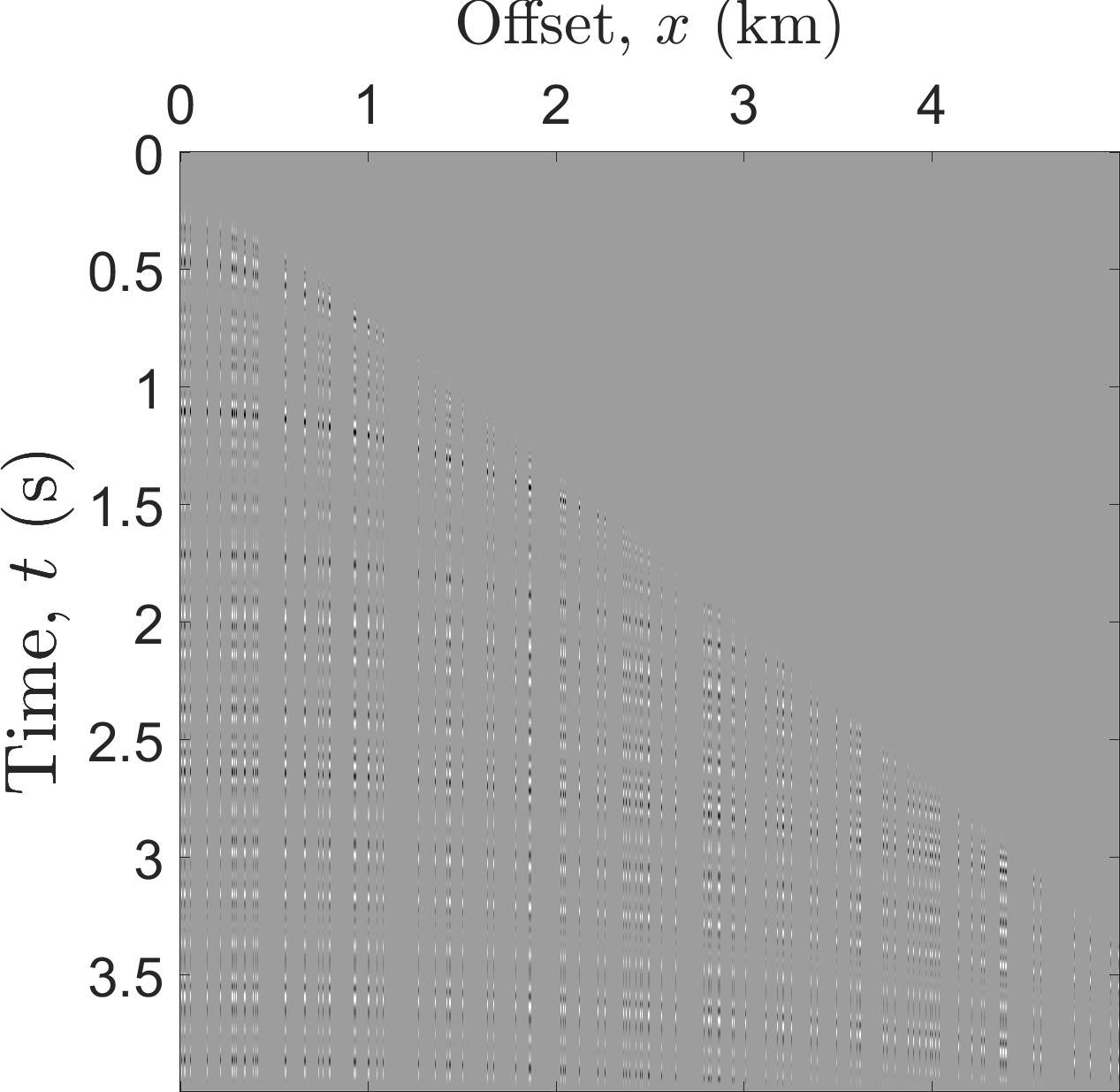}
		\caption{90\% missing traces}
	\end{subfigure}	
	\begin{subfigure}[b]{0.425\textwidth} \includegraphics[width=1\textwidth]{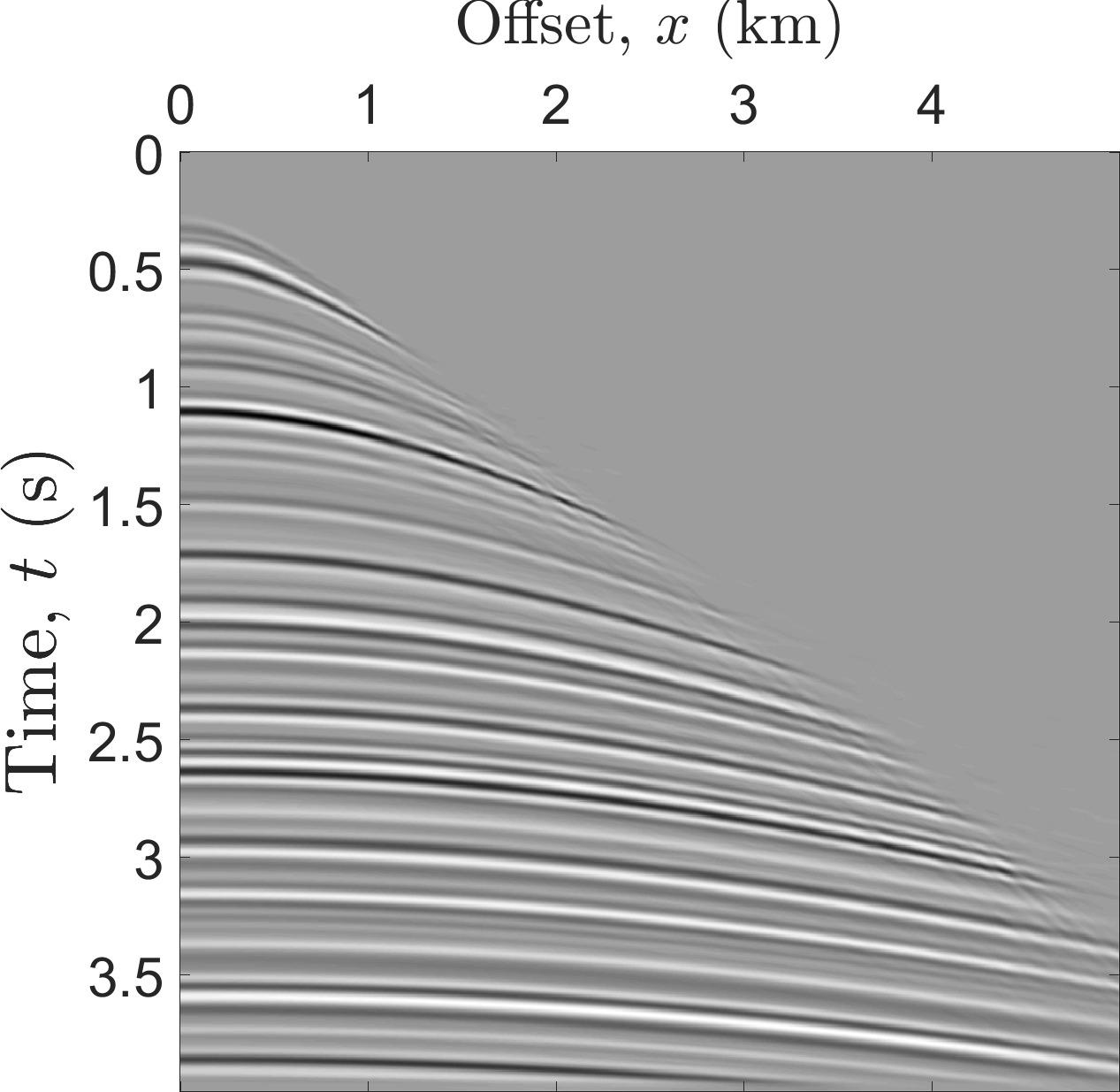}
		\caption{Reconstruction, $10\times$ thresholding}		
	\end{subfigure}			
	
	\caption{Interpolation into missing traces with the soft-thresholding algorithm.}
	\label{fig:cmpmissed2}
\end{figure}
Here we show some examples where we use soft-thresholding for conducting interpolation in cases of missing traces in the sampling setups. The CMP gather in Figure \ref{fig:cmpmissed}a contains 50\% randomly missing traces. For the data reconstruction we use the simple modification \eqref{itersmd} of the iterative scheme for obtaining sparse representations. Figure \ref{fig:cmpmissed}b shows reconstruction results after 30 soft-thresholding iterations. Note the absence of high amplitude artifacts produced by the proposed method. To control the obtaining results, we also consider synthetic CMP gathers with $90\%$ of missing traces, see figures \ref{fig:cmpmissed2}a,c. In spite of the low amount of given data, it is still possible to reconstruct the structure of the waves (Figure \ref{fig:cmpmissed2}b). Moreover, varying the parameter of soft-thresholding ($\mu$, see Chapter \ref{rec_tec}), one can improve the reconstruction quality. Here, the increase of the parameter $\mu$ leads to a better accuracy for low-amplitude events; conversely, high-amplitude events can be reconstructed with smaller values of $\mu$. In Figure \ref{fig:cmpmissed2}d we show the result of the reconstruction with soft-thresholding iterations where the parameter $\mu$ was increased by 10 times compared to the one used for reconstructions in Figures \ref{fig:cmpmissed}b and Figure \ref{fig:cmpmissed2}b.

\subsection{2D field data.}
\begin{figure}[t!]
	\centering 	
	\begin{subfigure}[b]{0.425\textwidth} \includegraphics[width=1\textwidth]{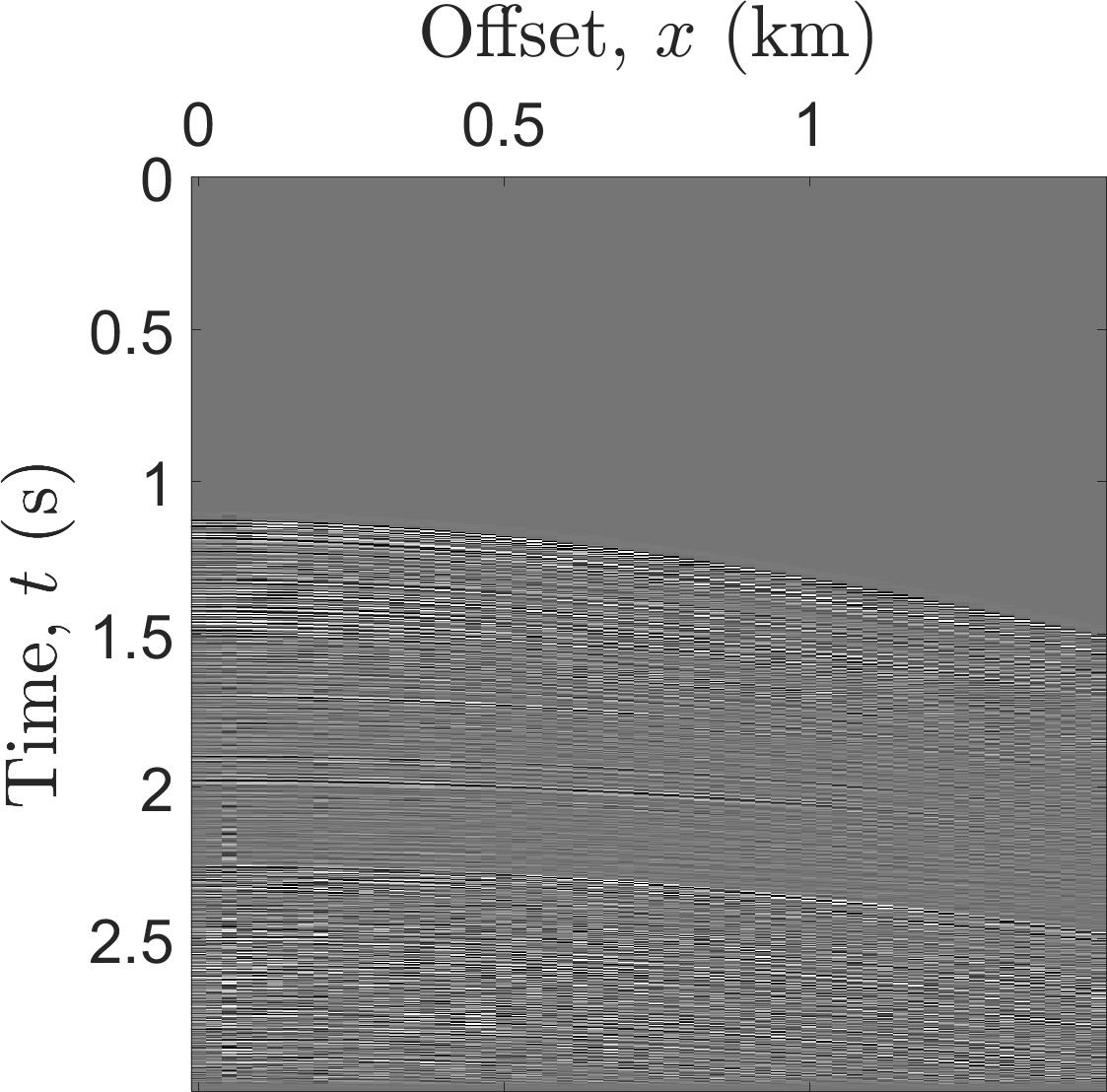}
		\caption{CMP gather}
	\end{subfigure}\hspace{5pt}	
	\begin{subfigure}[b]{0.425\textwidth} \includegraphics[width=1\textwidth]{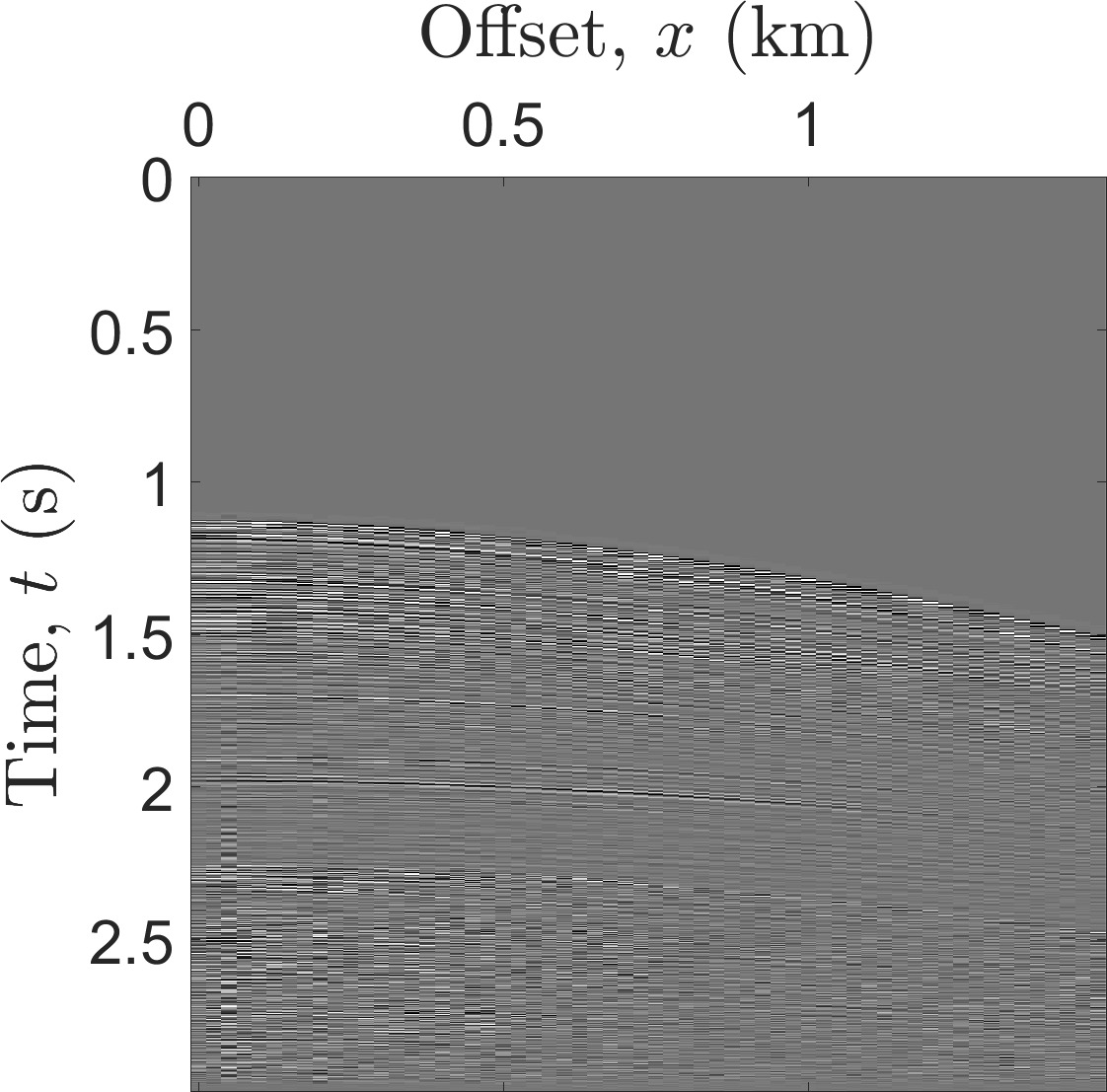}
		\caption{Subtracted multiples}
	\end{subfigure}	
	
	\vspace{5pt}	
	
	\begin{subfigure}[b]{0.425\textwidth} \includegraphics[width=1\textwidth]{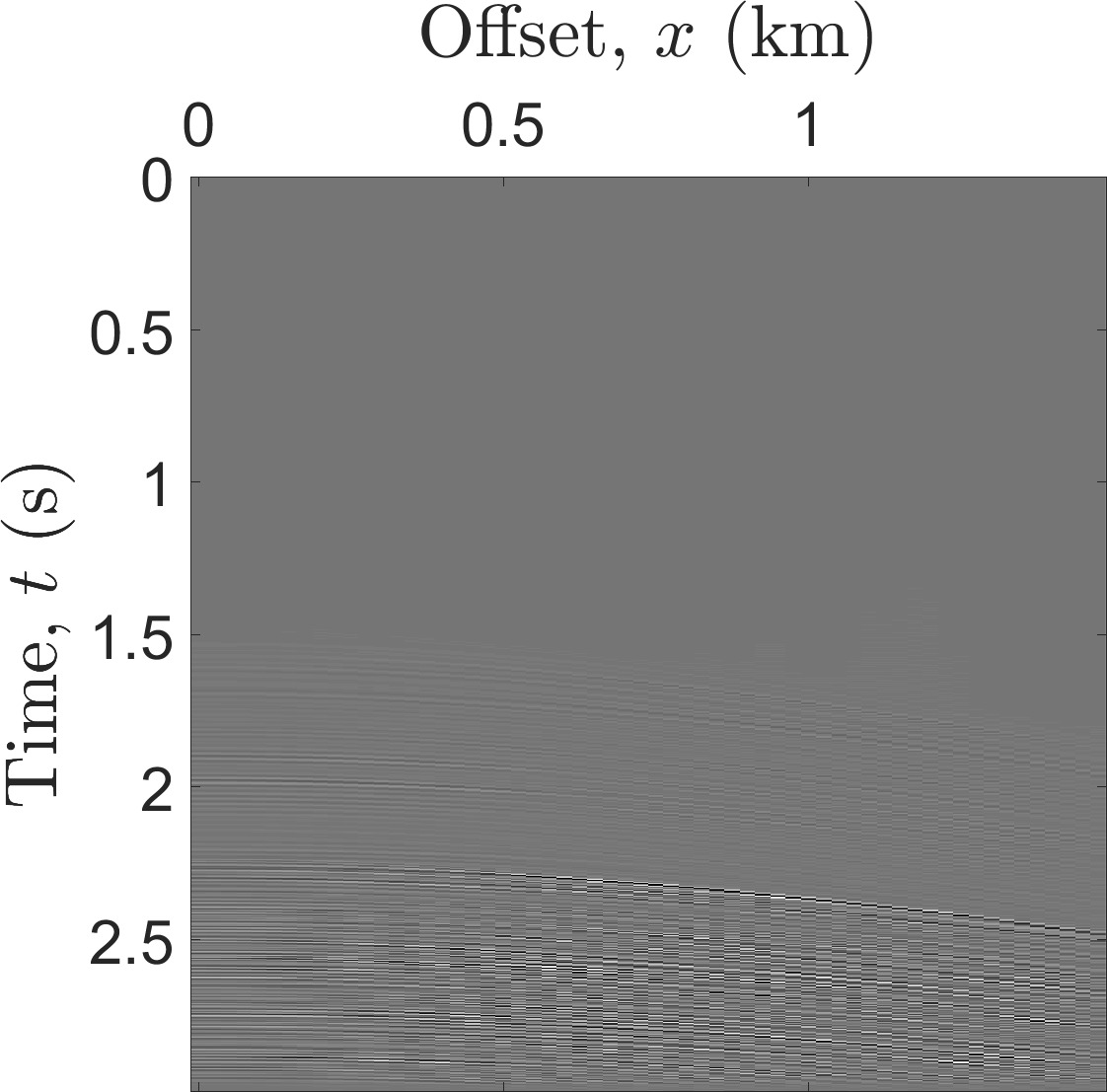}
		\caption{Reconstructed multiples after the sparse representation}
	\end{subfigure}\hspace{5pt}		
	\begin{subfigure}[b]{0.425\textwidth} \includegraphics[width=1\textwidth]{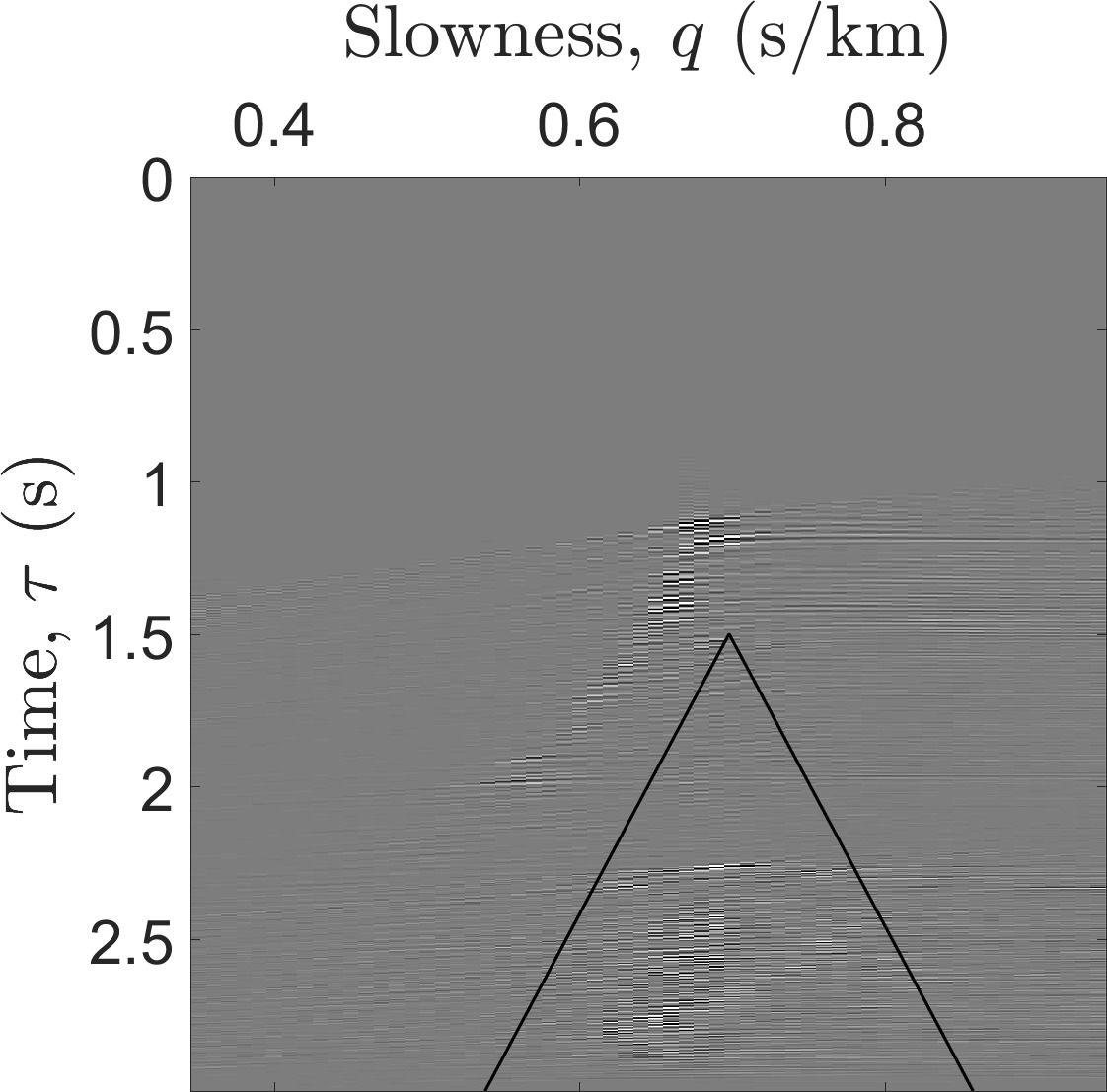}
		\caption{Result of soft-thresholding iterations, muting}
	\end{subfigure}		
	\caption{Multiple attenuation for 2D field data}
	\label{fig:cmpitersre}
\end{figure}
As an example of real data processing, we consider a CMP gather from the Canterbury data set \cite{lu2003three}. Multiple reflections in this CMP gather start at around 2.2 s (Figure \ref{fig:cmpitersre}a). Attenuation of the multiples was carried out after applying the reconstruction method from Chapter \ref{rec_tec} with 40  soft-thresholding iterations and the related muting procedure (Figure \ref{fig:cmpitersre}d). The part of the Radon image corresponding to multiples was taken back to the time-offset domain (Figure \ref{fig:cmpitersre}c) and subtracted from the initial CMP gather (Figure \ref{fig:cmpitersre}b).

\section{Conclusions}
A fast log-polar-based method for the evaluation of the hyperbolic Radon transform has been presented. According to the tests performed, the method demonstrates reasonable accuracy and favorable computational costs compared to other methods. The accuracy of the method can be increased when considering higher order interpolation kernels for coordinate conversions between time-offset, Radon, and log-polar domains. Numerical tests show that the GPU implementation is more than 10000 faster for large data sets in comparison to a direct implementation in standard C of sums over hyperbolas, and a substantial speedup is also obtained compared to alternative fast methods.

\clearpage
\section*{Acknowledgements}
This work was supported by the Crafoord Foundation (2014-0633) and the Swedish Research Council (2011-5589, 2015-03780)


\bibliographystyle{plain} 
\bibliography{hypradon_arxiv}


\end{document}